\documentclass[11pt, letterpaper]{amsart}
\usepackage{amsmath, amsthm, amscd, amsfonts, amssymb, graphicx, xcolor}
\usepackage[pagebackref=true]{hyperref}
\usepackage[alphabetic]{amsrefs}
\usepackage[english]{babel}
\usepackage[utf8]{inputenc}
\usepackage{lmodern}
\usepackage[T1]{fontenc}
\usepackage{enumitem}
\usepackage{stmaryrd}
\usepackage{floatrow}
\usepackage{dsfont}
\usepackage{subeqnarray}
\usepackage[all]{xy}
\usepackage{varioref}
\usepackage{tikz-cd}
\usetikzlibrary{arrows}

\usetikzlibrary{decorations.pathreplacing,decorations.markings}

\usepackage{pgfplots}\pgfplotsset{compat=1.18}
\usepackage[toc]{appendix}

\tikzset{
	on each segment/.style={
		decorate,
		decoration={
			show path construction,
			moveto code={},
			lineto code={
				\path [#1]
				(\tikzinputsegmentfirst) -- (\tikzinputsegmentlast);
			},
			curveto code={
				\path [#1] (\tikzinputsegmentfirst)
				.. controls
				(\tikzinputsegmentsupporta) and (\tikzinputsegmentsupportb)
				..
				(\tikzinputsegmentlast);
			},
			closepath code={
				\path [#1]
				(\tikzinputsegmentfirst) -- (\tikzinputsegmentlast);
			},
		},
	},
	mid arrow/.style={postaction={decorate,decoration={
				markings,
				mark=at position .5 with {\arrow[#1]{stealth}}
	}}},
}
\newcommand{\cat}{\text{CAT}}

\newcommand{\Map}{\operatorname{Map}}
\newcommand{\Mapstat}{\operatorname{Map}^{\operatorname{stat}}}
\newcommand{\bnd}{\operatorname{bnd}}
\newcommand{\bdd}{\operatorname{Bdd}}

\newcommand{\Xhu}{\widetilde{X}^\infty}
\newcommand{\Xh}{\widetilde{X}}
\newcommand{\Xhb}{\widetilde{X}^F}
\newcommand{\bdg}{\partial_{\text{Grom}}}
\newcommand{\Xhuc}{\overline{X}^\infty}
\newcommand{\Xhc}{\overline{X}}
\newcommand{\Xhbc}{\overline{X}^F}

\DeclareMathOperator{\supp}{\text{supp}}

\DeclareMathOperator{\iso}{Isom}
\newcommand{\bbR}{\mathbb{R}}
\newcommand{\bbP}{\mathbb{P}}
\newcommand{\bbZ}{\mathbb{Z}}
\newcommand{\bbN}{\mathbb{N}}

\newcommand{\calF}{\mathcal{F}}
\newcommand{\calB}{\mathcal{B}}

\newcommand{\proj}{\operatorname{proj}}

\newcommand{\Id}{\operatorname{Id}}
\newcommand{\bd}{\partial_\infty}
\newcommand{\nui}{\check{\nu}}
\newcommand{\Ti}{\check{T}}
\newcommand{\Bi}{\check{B}}
\newcommand{\Si}{\check{S}}
\newcommand{\Omegai}{\check{\Omega}}
\newcommand{\omegai}{\check{\omega}}
\newcommand{\alphai}{\check{\alpha}}
\newcommand{\bbPi}{\check{\mathbb{P}}}

\newcommand{\prob}{\operatorname{Prob}}
\newcommand{\invmeas}{\mathcal{I}_\mathbb{P}(\varphi)}

\newcommand{\probP}{\operatorname{Prob}_\mathbb{P}}
\newcommand{\XG}{\overline{X}^{\text{Grom}}}
\theoremstyle{plain}
\newtheorem{thm}{Theorem}[section]

\newtheorem{cor}[thm]{Corollary}
\newtheorem{lem}[thm]{Lemma}
\newtheorem{prop}[thm]{Proposition}
\theoremstyle{definition} 
\newtheorem{Def}[thm]{Definition}
\newtheorem{rem}[thm]{Remark}
\theoremstyle{definition} 

\newtheorem{ex}[thm]{Example}
\theoremstyle{definition}
\newtheorem*{ackn}{Acknowledgements}

\newcommand{\newcomment}[4]{%
	\newcounter{#2counter}
	\expandafter\newcommand\csname #1\endcsname[1]{%
		\refstepcounter{#2counter}%
		{\color{#4}(#3\arabic{#2counter})}\marginpar{\scriptsize\raggedright\textbf{\color{#4}(#2 \arabic{#2counter}):} ##1}%
}}

\newcomment{clb}{Corentin}{c}{teal}

\newcommand\extrafootertext[1]{%
	\bgroup
	\renewcommand\thefootnote{\fnsymbol{footnote}}%
	\renewcommand\thempfootnote{\fnsymbol{mpfootnote}}%
	\footnotetext[0]{#1}%
	\egroup
}
\title{Ergodic cocycles in hyperbolic and Hadamard spaces}
\author{Corentin Le Bars}
\begin{document}
	\begin{abstract}
		We consider discrete random dynamical systems induced by a non-elementary group action on a non-proper hyperbolic space. We prove that if the system is ergodic and satisfies the ``asymptotic past and future independence condition'' as defined by Bader and Furman \cite{bader_furman14}, the associated ergodic cocycle converges to the Gromov boundary almost surely. If the cocycle has finite first moment, we show that its drift is positive. Using hyperbolic models introduced in \cite{petyt_spriano_zalloum24}, we prove analogous statements for groups acting on Hadamard spaces with a pair of contracting elements. 
	\end{abstract}

	\maketitle
	\setcounter{tocdepth}{1}
	\tableofcontents
	\extrafootertext{The author is supported by ERC Advanced Grant NET 101141693.}
	\section{Introduction}
	
	Studying the asymptotic behavior of dynamical systems induced by group actions has been a very active area of research over the last 50 years, and has lead to major developments in group theory. This subject is related to a wide range of questions such as harmonic analysis, growth and amenability properties of $G$ and rigidity results. An emblematic example is the study of products of random matrices, after the seminal works of H.~Furstenberg and H.~Kesten \cites{furstenberg_kesten60,furstenberg63}. Another context of importance is that of hyperbolic groups and more generally, of groups acting on spaces of negative curvature. Random walks in these spaces typically escape to infinity along well-defined directions and satisfy strong limit laws, see \cite{kaimanovich00} for the convergence of random walks in hyperbolic groups or \cite{karlsson_margulis99} for $\cat$(0) spaces. The subject is still active and has known spectacular progress recently, see \cites{blachere_haissinsky_mathieu11,maher_tiozzo18,gouezel22,choi23} to name just a few. In this paper, we study random dynamical systems induced by a group acting on spaces sharing these hyperbolic properties.

	\subsection{Statement of the results}
	Let $(\Omega, \mathcal{F}, \mathbb{P})$ be a standard probability space, and $T : \Omega \to \Omega$ be a measure-preserving, invertible and ergodic transformation. Let $G$ be a discrete countable group acting by isometries on a metric space $(X,d)$. Let $f: \Omega \to G$ be a measurable map. 
	The random dynamical system generated by the data $(\Omega, \mathbb{P}, T, G, f)$ can be completely described by $\varphi(0, \omega) = \Id $ and for $n \geq 1, \omega \in \Omega$,
	\[\left\{ \begin{array}{ll}
		\varphi(n, \omega) = f(T^{n-1}\omega)\dots f(\omega) \nonumber \\
		\varphi(-n, \omega) = f(T^{-n}\omega)^{-1}f(T^{-n +1}\omega)^{-1}\dots f(T^{-1}\omega)^{-1} \nonumber.
	\end{array}\right.
	\]

	The cocycle $(\varphi(n, \omega))_{n\in \mathbb{Z}}$ describes a stochastic walk where the increments are distributed according to $f_\ast \mathbb{P}$, but are not necessarily independent. This situation therefore contains and extends the random walk setting. The associated random dynamical system on $X$ is 
	\begin{align}
		\varphi \colon& \mathbb{Z} \times \Omega \times X \to X \nonumber \\
		& (n, \omega,x) \mapsto \varphi(n, \omega)x,\nonumber 
	\end{align}
	with the induced $G$-action on $X$.
	
	We denote by $\chi$ the backward random cocycle (the term ``backward'' will be explained in Section \ref{section rds}), defined by 
	$$\chi(n, \omega) = \varphi(-n, \omega)^{-1} = \varphi(n, T^{-n}\omega).$$
	In particular, for all $n \in \mathbb{N}$, 
	$$ \chi(n, \omega) = f(T^{-1}\omega)f(T^{-2}\omega)\dots f(T^{-n}\omega).$$
	In the random walk case where $\mu \in \prob(G)$ is a probability measure on $G$, the space $(\Omega = G^\mathbb{Z}, \mathbb{P}= \mu^{\otimes \mathbb{Z}})$ is the space of increments and $T$ is the Bernoulli shift, the cocycle $\chi$ is the classical right random work product of independent, identically $\mu$-distributed elements in $G$:
	$$\chi(n, \omega) = \omega_{-1}\omega_{-2}\dots \omega_{-n}.$$
	
	A key feature of Furstenberg's approach to the study of random walks is the use of boundary theory, i.e. the study of measures supported on well-chosen boundaries and which are ``invariant'' under the random dynamical system. The precise notion for this invariance in the context of random walks is that of \emph{stationarity}.
	We adopt this approach, although, ergodic cocycles being time-dependent, the notion of invariance has to be adapted. As Furstenberg's approach relies on martingale convergence theorems in a crucial way, we emphasize that it is natural to consider products on the right and the backward cocycle $\chi$. 
	
	In this paper, we consider the two following situations: 
	\begin{itemize}
		\item $X$ is a (possibly non-proper) almost-geodesic Gromov-hyperbolic space, and a discrete countable group $G$ acts on $X$ with two independent loxodromic elements, i.e. loxodromic elements with disjoint fixed points at infinity;
		\item $X$ is a (possibly non-proper) complete $\cat$(0) space, and $G$ acts on $X$ with two independent contracting elements. 
	\end{itemize}
	
	Loxodromic elements in hyperbolic spaces and contracting elements in $\cat$(0) spaces are elements with a strong hyperbolic behavior: for instance, they induce a North-South dynamics on respectively the Gromov boundary $\bdg X$ or the visual boundary $\bd X$. In the rest of the paper, we say that the $G$-action on $X$ (for $X$ either Gromov-hyperbolic or $\cat$(0)) is \emph{non-elementary} when $G$ possesses such a pair of elements for this action. Note that in the case of an action on a proper $\cat$(0) space, contracting elements are exactly rank one elements \cite{bestvina_fujiwara09}. Moreover, the Rank Rigidity Conjecture states that ``generic'' geometric actions (that is, where $X$ is not a higher rank symmetric space nor a building of dimension $\geq 2$) on irreducible geodesically complete proper $\cat(0)$ spaces contain rank one elements. 
	
	Obviously, in order to use the non-elementarity of the action and to prove some results of convergence, we need to add some conditions on the dynamical system (the function $f$ could be essentially constant). The condition we use here is that of \emph{asymptotic past and future independence}, which roughly says that the behavior of the random dynamical system in the asymptotic future is independent from its behavior in the asymptotic past. This property can be thought of as a (measurable) \emph{local product structure} for the dynamical system, similar to what happens for the geodesic flow on compact hyperbolic surfaces (or more generally, Anosov flows): given a past trajectory and a future trajectory of the geodesic flow, there exists a geodesic line which has these prescribed past and future asymptotically. The asymptotic past and future condition (we say that the random dynamical system is ``apafi'') was introduced by U.~Bader and A.~Furman in \cite{bader_furman14}, and is central to our argument. 
	
	For the rest of this section, we let $(\Omega, \mathcal{F}, \mathbb{P}, T,f, G)$ be a discrete random dynamical system on a discrete countable group $G$ such that the measure-preserving flow $T: \Omega \to \Omega$ is ergodic. We make the assumption that the cocycle $\varphi : \mathbb{Z} \times \Omega \rightarrow G$ satisfies the asymptotic past and future independence condition. 
	The first result is the following convergence theorem in hyperbolic spaces. 
	
	\begin{thm}[{Convergence in hyperbolic spaces}]\label{thm cv cocycle hyp intro}
		Let $(X,d)$ be a separable almost-geodesic Gromov-hyperbolic space and $G \curvearrowright (X, d)$ be a non-elementary action by isometries. Then, for any $o\in X$, the backward cocycle $(\chi(n, \omega) o)_n$ converges $\mathbb{P}$-almost surely to a point $\xi^+(\omega) \in \bdg X$ in the Gromov boundary. 
	\end{thm}
	
	Notice that the convergence takes place without integrability assumption on the cocycle $\varphi$. Notice also that we do not require $X$ to be proper, so this result applies to the following situations: 
	\begin{itemize}
		\item an isometric action on a $\cat(-1)$ space, for instance the Cremona group on the Picard-Manin space;
		\item the action of the mapping class group $\operatorname{Mod}(S_{g,p})$ of a surface on its arc graph and curve complex; 
		\item actions on $\bbR$-trees, in particular trees associated with rank-1 valuations on fields;
		\item the action of $\operatorname{Out}(F_N)$ on the free splitting complex $\mathcal{FS}(F_N)$ and the free factor complex $\mathcal{FF}(F_N)$;
		\item the action of a right-angled Artin group on its extension graph.
	\end{itemize}
	More generally, any acylindrically hyperbolic group satisfies this condition. If the RDS is in fact a $\mu$-generated random walk on $G$ (for $\mu$ admissible), we recover the convergence result of Maher and Tiozzo \cite{maher_tiozzo18}. Since we only consider actions of discrete countable groups, the assumption of separability can also be dropped, see Section \ref{section separable} below. 
	
	We recall that a Hadamard space is a metrically complete $\cat$(0) space. The analogue of Theorem \ref{thm cv cocycle hyp intro} in the context of Hadamard spaces is the following.

	\begin{thm}[{Convergence in Hadamard spaces}]\label{thm cv cocycle cat intro}
		Let $(X,d)$ be a separable Hadamard space and $G \curvearrowright (X, d)$ be an action by isometries and assume that $G$ contains a pair of independent contracting isometries for this action. Then, for any $o\in X$, the backward cocycle $(\chi(n, \omega) o)_n$ converges $\mathbb{P}$-almost surely to a point $\xi^+(\omega) \in \bd X$ in the visual boundary. 
	\end{thm}
	
	The proof of Theorem \ref{thm cv cocycle cat intro} relies on the use of well-chosen hyperbolic models $X_L:=(X, d_L)_{L \in \mathbb{N}}$ attached to the Hadamard space $(X,d)$. These models are due to H.~Petyt, D.~Spriano and A.~Zalloum \cite{petyt_spriano_zalloum24}, and rely on the combinatorial structure induced by \emph{curtains}, which can be thought of as hyperplanes of $X$. This construction is similar to that of the contact graph for $\cat$(0) cube complexes \cite{hagen14}, and is inspired by ideas of Genevois \cite{genevois19}. In the case of the contact graph associated to a $\cat$(0) cube complex $Y$, there exists an $\iso(Y)$-equivariant embedding of the Gromov boundary of the contact graph into the Roller boundary $\partial_{\mathcal{R}}Y$ \cite[Theorem 3.7]{fernos_lecureux_matheus21} of $Y$. We prove the analogue in the context of general Hadamard spaces, which was known for locally compact $\cat$(0) spaces \cite[Theorem~8.1]{petyt_spriano_zalloum24}. 
	
	\begin{thm}\label{thm homeomorphism bords XL intro}
		Let $X$ be a complete $\cat$(0) space. Then for any $L\geq 0 $, there is an $\iso(X)$-equivariant homeomorphic embedding of the boundary $\partial_L : \bdg X_L \rightarrow \bd X$. 
	\end{thm}

	Thanks to this homeomorphism, Theorem \ref{thm cv cocycle cat intro} is just a consequence of Theorem \ref{thm cv cocycle hyp intro} when applied to a well-chosen hyperbolic model.

	We can then derive some non-commutative laws of large numbers. Given a RDS on $X$ as before, we say that the cocycle is integrable if for some (equivalently, any) $o \in X$,
	$$\int_\Omega d(\varphi(1, \omega)o, o) d\mathbb{P}(\omega) < \infty.$$ 
	Let $\varphi$ be such an integrable cocycle on $X$, and consider the displacement variables 
	$(D(n, \omega))_n$, where
	$$ D(n, \omega) := d(\chi(n, \omega)o , o), $$ 
	and define the asymptotic average 
	$$\lambda = \lambda_X(\varphi) := \inf \frac{1}{n} \int_\Omega D(n, \omega) d\mathbb{P}(\omega) \in [0, +\infty).$$
	By a direct application of the Kingman subadditive Theorem \cite{kingman68}, the variables $(D(n, \omega)/n)_n$ converge almost surely to $\lambda$ and the convergence also holds in $L^1(\mathbb{P})$. 
	
	\begin{thm}[{Positive drift in hyperbolic spaces}]\label{thm drift hyp intro}
		Let $(X,d)$ be a separable almost-geodesic Gromov-hyperbolic space and $G \curvearrowright (X, d)$ be a non-elementary action by isometries. Assume that the cocycle $\varphi$ is integrable. Then $\lambda_X(\varphi) > 0$. 
	\end{thm}
	
	Thanks to the hyperbolic models, the analogue holds for Hadamard spaces. 
	
	\begin{thm}[{Positive drift in Hadamard spaces}]\label{thm drift cat intro}
		Let $(X,d)$ be a separable Hadamard space and $G \curvearrowright (X, d)$ be an action by isometries admitting a pair of contracting elements. Assume that the cocycle $\varphi$ is integrable. Then $\lambda_X(\varphi) > 0$. 
		For $\mathbb{P}$-almost every $\omega \in \Omega$, there is a unique geodesic ray $\gamma^\omega : [0, \infty) \rightarrow X$ starting at $o$ such that 
		\begin{equation}
			\lim_{n\rightarrow \infty} \frac{1}{n} d(\gamma^\omega(\lambda_X n), \chi(n, \omega)o) = 0. \nonumber
		\end{equation} 
	\end{thm}
	
	A series of recent works aim at studying groups acting with contracting elements on metric spaces in the greatest generality. Paradigmatic questions in this context are: is there a good notion of boundary and can we endow this boundary with Patterson-Sullivan measures \cite{coulon24,yang25}; are there isometry-equivariant hyperbolic models on which contracting elements in $G \curvearrowright X$ become loxodromic elements, and what are the geodesic tracking properties of such models \cite{petyt_zalloum24,zbinden24}. The approach of limit laws using boundary theory developed in this paper is particularly well-adapted to this theme. In particular, the papers \cite{petyt_zalloum24,zbinden24} prove that in a very general context (that encompasses many acylindrically hyperbolic actions, injective spaces, Garside groups, HHS...), one can construct such a hyperbolic model with good recognition of contracting elements. Therefore Theorems \ref{thm cv cocycle hyp intro} and \ref{thm drift hyp intro} apply. In order to translate this information back to the original space $X$, one needs to study how the hyperbolic directions in the model correspond to points of a well-chosen boundary, for instance the sublinear Morse boundary, as in Theorem \ref{thm homeomorphism bords XL intro}. Once we have such a tracking, the approach presented here and the limit laws would extend to this more general context. We plan to treat this question in a future work. 
	\newline 
	
	Along the way, we are lead to consider some auxiliary measurable spaces associated to the  dynamical system: the spaces of ideal futures and of ideal pasts. Each of them can be defined as the Mackey range of a suitably chosen cocycle, and they generalize the past and future Poisson-Furstenberg boundaries for random walks (see Section \ref{section space trajectories} for the definition). We derive several properties that are reminiscent of classical results for Poisson-Furstenberg boundaries, and show that they satisfy strong ergodic properties. In particular, we show that under the apafi condition, the $G$-action on these spaces is coarsely metrically ergodic, as defined in \cite{bader_caprace_furman_sisto22}.

	\begin{thm}[{Coarse isometric ergodicity of the Mackey range}]\label{thm coarse isom erg intro}
		Let $G$ be an locally compact second countable group and let $(\Omega, \mathcal{F}, \mathbb{P}, T, f, G)$ be an ergodic invertible RDS  on $G$ satisfying the apafi condition. Let $B_-$ and $B_+$ be the corresponding past and future Mackey ranges associated to $\varphi$. Then the diagonal $G$-action on $B_- \times B_+$ is metrically ergodic and coarsely metrically ergodic. 
	\end{thm}
	We refer to Section \ref{section erg prop mackey range} for the definition of coarse metric ergodicity. Although metric ergodicity of the $G$-action on $B_- \times B_+$ was known \cite{bader_furman14}, the novelty here is the fact that the action is also coarsely metrically ergodic. 
	
	In order to prove this theorem, we are lead to show that the action on any of these spaces is \emph{strongly approximately transitive} as defined by Jaworski \cite{jaworski95}. This notion was inspired by the proximal properties of random walks, and applied to give an alternative proof of the fact that virtually nilpotent groups are Choquet-Deny, i.e. do not have non-constant bounded harmonic functions. The Poisson-Furstenberg boundary $B(G, \mu))$ of a probability measure $\mu$ on a discrete countable group $G$ is SAT \cite[Corollary 2.4]{jaworski95}, and we show the analogous statement in the context of ergodic cocycles. 
	
	\begin{prop}\label{prop sat intro}
		The $G$-action on $(B_\pm, \nu_\pm)$ is strongly approximately transitive. 
	\end{prop}
	For a $\mu$-generated random walk, the Poisson boundary $(B_\mu, [\nu_B])$ is characterized by the Poisson formula: there is a surjective isometry 
	\[\Phi_\mu\colon L^\infty(B_\mu, [\nu_B]) \to \operatorname{Har}^\infty(G, \mu)\] 
	from the set of essentially bounded functions on $B_\mu$ to the set of bounded $\mu$-harmonic functions on $G$. 
	The key tool in the proof of Proposition \ref{prop sat intro} and Theorem \ref{thm coarse isom erg intro} is a version of the Poisson transform  adapted to our context. 
	
	\begin{ex}[Compact negatively curved manifolds]
		Let $M$ be a compact hyperbolic manifold, and $\tilde{M}$ be its universal covering space. There is a one-to-one correspondence between invariant measures $\bbP$ of the geodesic flow on the unit tangent bundle  $\Omega = UM$ and \emph{geodesic currents}, i.e. Radon measures on the space $\partial \tilde{M} \times \partial \tilde{M} -\Delta(\partial \tilde{M} )$ (where $\Delta(\partial \tilde{M} )$ is the diagonal) which are invariant with respect to the action of the fundamental group $\pi_1(M)$, see \cite{kaimanovich94}. If $\bbP$ is a Gibbs measure, then the associated geodesic current is a \emph{quasi-product}, meaning that it is locally equivalent to the product of two $G$-invariant measure classes $[\bbP_-]$ and $[\bbP_+]$  on $\partial \tilde M$ that are given by the conditional measures of $\bbP$ on the strongly stable and strongly unstable foliations and of the Lebesgue measure on the geodesics. For instance, the maximal entropy (Bowen-Margulis) measure or the Lebesgue-Liouville measures are Gibbs measures. Ledrappier shows in \cite{ledrappier88} that the harmonic class (hitting measure of a Brownian motion on the covering space) is also Gibbs. Let now $T$ be the time-1 geodesic flow on $UM$, $\tilde{T}$ be the time-1 geodesic flow of the universal cover and $\sigma : UM \to U \tilde{M}$ be a measurable section of the covering map $\pi : U\tilde{M} \to U M$. Then the measurable map $f : \Omega \to \pi_1 (M)$ defined by 
		$$ \tilde{T}(\sigma(\omega)) = f(\omega) \sigma (T\omega)$$ 
		defines a RDS. The spaces of ideal future and ideal past $B_+$ and $B_-$ are realized on $\partial \tilde{M}$ and one has the apafi condition by the properties of the Gibbs measure that we chose \cite[Theorem~8.1]{bader_furman25}. 
		With this example in mind, it is good to compare the result of Proposition \ref{prop sat intro} with \cite{kaimanovich02}. 
	\end{ex}

	\begin{rem}
		Some parts of this paper were written prior to the publication of \cite{bader_furman25}. We kept the proofs that are particularly enlightening or paradigmatic, and for self-containment purpose, but we refer the reader to the aforementioned paper for a more detailed and general exposition to ergodic cocycles via boundary theory. 
	\end{rem}
	
	\subsection{Strategy}
	
	The strategy relies on the study of boundary maps induced by the spaces of ideal pasts and futures, which correspond to Poisson-Furstenberg boundaries in the random walk setting. Let us explain the case where $X$ is a hyperbolic space. In a similar way to the random walk case, the first step is to find a ``stationary'' measure: here it means a probability measure $\mu^+_\cdot \in \prob(\Omega_+ \times \bdg X)$ that is invariant under the forward random dynamical system. By a martingale argument, one obtains the existence of a family of probability measures $\{\mu_\omega\}_{\omega \in \Omega} \in \prob(\bdg X)$ such that we have the almost sure convergence 
	\begin{align}
		\chi(n, \omega) \mu^+_{T^{-n}\omega} \underset{n \to \infty}{\longrightarrow} \mu_\omega. \label{eq contraction intro}
	\end{align}
	In the same way as in the random walk case, these family of measures are shown to be Dirac measures $\delta_{\psi(\omega)}$, for $\psi: \Omega \to \bdg X$ a measurable map. This fact is a byproduct of the study of $G$-equivariant measurable maps $B_+ \to \prob (\bdg X)$, where $B_+$ is the space of ideal futures. Thanks to the strong ergodic properties of these boundaries, and notably the fact the the action is coarsely metrically ergodic as given by Theorem \ref{thm coarse isom erg intro}, we can appeal to the results of \cite{bader_caprace_furman_sisto22} that classify such $G$-maps. The contraction phenomenon given by Equation~\eqref{eq contraction intro} on the boundary can be turned into a statement of convergence in the space after a careful study of the support of the measures $\{\mu^{+}_{\omega}\}$. This contraction in probability corresponds of the boundary allows to conclude of the convergence thanks to weak convexity of shadows in hyperbolic spaces.

	\subsection{Structure of the paper}
	In Section \ref{section prelim hyp}, we give some preliminaries on hyperbolic spaces. As the hyperbolic models $\{X_L\}$ are not geodesic, we show that the results in \cite{bader_caprace_furman_sisto22}, which are central to our argument, remain valid. 
	In Section \ref{section hyperbolic models}, we review the construction of hyperbolic models for Hadamard spaces of Petyt-Spriano-Zalloum \cite{petyt_spriano_zalloum24} and we give the proof of Theorem \ref{thm homeomorphism bords XL intro} establishing a homeomorphism between the Gromov boundary of a hyperbolic model and a subset of the visual boundary of the $\cat$(0) space. 
	Section \ref{section rds} introduces discrete random dynamical systems and some important notions: disintegration of measures, invariant measures and restrictions to positive times. 
	We develop the theory of Mackey ranges of ergodic cocycles after Bader-Furman in Section \ref{section bd system}, and we show a number of their interesting ergodic properties. 
	In Section \ref{section erg coc hyp cv}, we collect the previous results to show that, under the apafi condition, the backward cocycle in a Gromov-hyperbolic space converges to a point of the Gromov boundary. 
	Section  \ref{section erg coc cat cv} is dedicated to the proof of the analogous result for Hadamard spaces. 
	Finally, Section \ref{section drift csq} studies the escape rate of the RDS: we prove that it is positive under an integrability assumption, and we derive a sublinear tracking result. 
	
	\begin{ackn}
		I am grateful to Uri Bader and Alex Furman for sharing their notes prior to their publication, and especially to Uri Bader for discussions. I am also thankful to the Weizmann Institute of Science for its support. Parts of this paper were written while I was visiting the Max Planck Institute for Mathematics in the Science, and I want to thank Anna Wienhard and the people there for their hospitality. 
	\end{ackn}
	
	\section{Preliminaries on hyperbolic spaces}\label{section prelim hyp}
	In this section we recall some basic facts about Gromov-hyperbolic spaces. What follows is standard, but since most papers deal about geodesic and proper hyperbolic spaces, we will present this more general setting. Standard references include \cite{bridson_haefliger99,drutu_kapovich18,bonk_schramm00,vaisala05}, the last three dealing also with non-geodesic spaces. 
	
	\subsection{First Definitions} 
	
	Let $(X,d)$ be a metric space. Given $x,y,o \in X$, the Gromov product between $x$ and $y$ with respect to $o $ is given by 
	\[(x|y)_o = \frac{1}{2} (d(w, o ) + d(y, o ) - d(x,y)).\]
	Note that if we take another basepoint $o' $, we have 
	\[|(x|y)_o - (x|y)_{o'} | \leq d(o, o').\]
	
	The following definition is due to Gromov \cite{gromov87}. 
	\begin{Def}
		We say that the metric space $(X,d)$ is \emph{$\delta$-hyperbolic} for $ 0\leq \delta < \infty$, if for all $x,y,z, o \in X$, 
		\[(x|y)_o \geq \min\{(x|z)_o, (z|y)_o\} - \delta.\]
		If $X$ is $\delta$-hyperbolic for some $\delta$, we say that it is \emph{Gromov-hyperbolic}. 
	\end{Def}
	
	Let $(X,d)$ be a Gromov-hyperbolic space, and let $o \in X$. We say that a sequence of points $(x_i)_{i \in \mathbb{N}}$ \emph{converges at infinity } if 
	\[\lim_{i,j \rightarrow \infty} (x_i|x_j)_o =  \infty,\]
	and that two sequences that converge at infinity are \emph{equivalent} if 
	\[\lim_{i \rightarrow \infty} (x_i|y_i)_o =  \infty.\]
	
	This is an equivalence relation, and as 
	$$|(x|y)_o - (x|y)_{0'}| \leq d(o, o'), $$
	this definition does not depend on the basepoint $o$. It is worth noticing that if $(x_n)$ and $(y_n)$ are two sequences that converge at infinity that are not equivalent, then the set $((x_i|y_j)_o)_{i,j}$ is bounded \cite[Lemma 5.3.(4)]{vaisala05}. 
	The Gromov boundary $\bdg X$ of $X$ is then defined as the set of equivalence classes of sequences $(x_i)$ that converge at infinity, for this equivalence relation. We denote by $\XG= X \cup \bdg X$ the \emph{Gromov bordification}, although for the moment no topology has been defined on it.

	We extend the Gromov product to points of the boundary in the following way. For $\xi, \eta \in \bdg X$, 
	\[(\xi| \eta)_o := \sup \{ \liminf_{i\rightarrow \infty} (x_i| y_i)_o \mid (x_i )\in \xi, (y_i )\in \eta\}.\]
	
	Taking the supremum and $\liminf$ in this definition is not the only choice. However, all four possible choices only differ by $2 \delta$ \cite[Lemma 5.6]{vaisala05}. Note that $(\xi|\eta)_o = \infty$ if and only if $\xi = \eta \in \bdg X$. Accordingly, if $x \in X, \xi \in \bdg X$, we let 
	\[	(\xi| x)_o := \sup \{ \liminf_{i\rightarrow \infty} (x| y_i)_o \mid (y_i ) \in \xi\}. \]
	
	As a consequence of the definitions, if $(X,d)$ is a $\delta$-hyperbolic space, there exists a constant $C= C(\delta)$ depending only on $\delta$ such that for all $x, y , z \in X \cup \bdg X$,
	\begin{align}
		(x|y)_o \geq \min\{(x|z)_o, (z|y)_o\} - C. \nonumber
	\end{align}

	The topology on $\bdg X $ is given by the basis of open sets 
	\begin{align}
		U_{\xi, R} = \{ x \in \XG \mid (\xi| x)_o > R\} \nonumber, 
	\end{align}
	so that $x_n \rightarrow \xi \in \bdg X$ if and only if $(\xi|x_n)_o \rightarrow \infty$. Moreover, if $(x_n)$ is a sequence in $X$, $(x_n)$ converges to $\xi \in \bdg X$ in this sense if and only if $(x_n) $ converges at infinity and the equivalence class $[(x_i)]$ of $(x_i)$ is $ \xi$ \cite[Lemma 11.101]{drutu_kapovich18}. For $Y\subseteq X$ a subset, we will write $\overline{Y}^{\text{Grom}}$ the closure of $Y$ in $\XG$ with this topology. The bordification  $X \cup \bdg X$ is completely metrizable and $\bdg X$ can be endowed with a family of metrics $d_\varepsilon$ depending on good choices of $\varepsilon >0$. These can be defined as 
	$$d_\varepsilon (\xi, \eta) = \inf \sum_{j=1}^n \exp (-\varepsilon (x_i| x_{i+1})_o),$$
	where the infimum is taken over all sequences $\xi = x_1, \dots , x_{n+1}= \eta$. This metric is compatible with the topology on $\XG$, and is complete \cite[Proposition 5.13]{vaisala05}. 
	
	\subsection{Coarse geometry}
	The hyperbolic spaces we are going to consider are not geodesic, so let's introduce some terminology in coarse geometry. 
	
	Let $f : X \rightarrow Y$ be a map between metric spaces, and let $Q \geq 1, C \geq 0$. We say that $f$ is a \emph{$(Q,C)$-quasi-isometric embedding} if for all $x, y \in X$, 
	\begin{align}
		\frac{1}{Q} d_X(x, y) - C \leq d_Y(f(x) , f(y)) \leq Q d_X(x,y) + C \nonumber. 
	\end{align}
	If moreover, there exists $A>0$ such that every point $y \in Y$ is at distance less than $A$ from $f(X)$, we say that $f$ is a \emph{$(Q,C)$-quasi-isometry}. In the particular case where $I \subseteq \bbR$ is an interval, and $f : I \rightarrow (X,d)$ is a $(Q,C)$-quasi-isometric embedding, we say that $f$ (or $f(I)$) is a \emph{$(Q,C)$-quasigeodesic path}. If $I = [a, b]$ is a compact interval, resp. $I = [a, \infty)$, resp. $I = \bbR$, we say that $f$ is a \emph{quasigeodesic segment}, resp. \emph{quasigeodesic ray}, resp. \emph{quasigeodesic line}. 
	
	Last, we say that a space is \emph{$(Q,C)$-quasigeodesic} if for every $x,y \in X$, there is a $(Q,C)$-quasigeodesic segment $\gamma : [a, b] \rightarrow X$ that joins $x$ to $y$, i.e. such that $\gamma(a) = x$ and $\gamma(b)= y  $. 
	
	\begin{Def}
		A metric space $(X,d)$ is \emph{$\alpha$-almost geodesic} if for every $ x, y \in X$, and every $t \in [0, d(x,y)]$, there is $z \in X $ such that 
		\begin{align}
			|d(x,z) - t | \leq \alpha \text{ and } |d(y, z ) - (d(x,y) - t)|\leq \alpha \nonumber.
		\end{align}
	\end{Def}
	
	In almost geodesic hyperbolic spaces, we have the following. 
	
	\begin{prop}[{\cite[Proposition 5.2]{bonk_schramm00}}]\label{prop exist qgeod hyp}
		Let $(X,d)$ be a $\delta$-hyperbolic, $\alpha$-almost geodesic metric space. Then there exists a universal constant $C= C(\delta, \alpha)$ such that:
		\begin{enumerate}
			\item \label{rough geodesic} for all $ x, y \in X$, there exists a $(1, C)$-quasigeodesic segment $\gamma : [a, b] \rightarrow X$ such that $\gamma(a) = x$ and $\gamma(b) = y$. 
			\item for all $ x\in X, y \in \bdg X$, there exists a $(1, C)$-quasigeodesic ray $\gamma : [0, \infty) \rightarrow X$ such that $\gamma(0) = x$ and $\lim_{t \rightarrow \infty} \gamma(t) = y$. 
			\item for all $ x, y\in \bdg X$, $x\neq y$, there exists a $(1, C)$-quasigeodesic line $\gamma : \bbR \rightarrow X$ such that $\lim_{t \rightarrow -\infty} \gamma(t)= x$ and $\lim_{t \rightarrow \infty} \gamma(t) = y$. 
		\end{enumerate}
	\end{prop}
	
	A metric space satisfying the conclusion (\ref{rough geodesic}) of the previous proposition is said to be \emph{$C$-roughly geodesic}. The spaces we are going to deal with are $\delta$-hyperbolic and $\alpha$-almost geodesic, for some $\delta, \alpha \geq 0$. In particular, they will be roughly geodesic, hence quasigeodesic. 
	
	In a geodesic $\delta$-hyperbolic space, there is a universal constant $C' = C'(\delta, Q, C)$ such that if $\gamma$ is a $(Q,C)$-quasigeodesic segment, then any geodesic segment between $x = \gamma(a)$ and $y = \gamma(b) $ lies in the $C'$-neighbourhood of $\gamma$. This phenomenon is called \emph{geodesic stability}. This result remains true for general hyperbolic spaces due to the following theorem of Bonk and Schramm. 
	
	\begin{thm}[{\cite[Theorem 4.1]{bonk_schramm00}}]\label{thm isom embed geod hyp space}
		Let $(X,d)$ be a $\delta$-hyperbolic space. Then there is an isometric embedding $ \iota : X \rightarrow Y $ of $X$ into a complete $\delta$-hyperbolic geodesic metric space $Y$.
	\end{thm}
	Moreover, $\iso(X)$ acts on the geodesic $\delta$-hyperbolic space by isometries and the embedding is $\iso(X)$-equivariant \cite[Corollary~A.10]{blachere_haissinsky_mathieu11}. 
	
	The analogue of geodesic stability in the general context is then the following, where we denote by $d_\text{Haus} (U, V) $ the Hausdorff distance between two subsets of $(X, d) $. 
	
	\begin{cor}[{\cite[Proposition 5.4]{bonk_schramm00}}]\label{cor stab qgeod hyp}
		Let $(X,d)$ be a $\delta$-hyperbolic metric space and $\gamma_1, \gamma_2$ be $(Q,C)$-quasigeodesic paths in $X$ with the same endpoints. Then there exists $C' = C'(\delta, Q, C)$ such that $d_\text{Haus}(\gamma_1, \gamma_2) \leq C'$. 
	\end{cor}

	We define the \emph{$R$-shadow} of a point $x $ seen from $o \in X$ by 
	\begin{align}
		S_o (x, R) := \{ y \in X \; : \; (y|o)_x \leq R\}. \nonumber
	\end{align}
	Rewriting the inequality of the Gromov product gives that for all $y, z \in \overline{S_o (x, 2R)}^{\operatorname{Grom}}$,
	\begin{align}
		(y|z)_o \geq d(o, x) - R + O(\delta) \nonumber, 
	\end{align}
	where $O(\delta)$ is a universal constant depending only on $\delta$. 
	
	A very useful fact about shadows in hyperbolic metric spaces is that they are weakly convex in the following sense.

	\begin{prop}\label{cor conv shadows}
		Let $(X,d)$ be a $\delta$-hyperbolic space. Then for every $Q\geq 1, D\geq 0$, there exists a constant $D = D(\delta, Q, C)$ such that:
		\begin{enumerate}
			\item if $y, z$ belong to the shadow $S_o (x, R)$, then any $(Q,C)$-quasigeodesic segment between $y$ and $z$ belongs to $S_o (x, R+D)$. 
			\item if $y, z$  do not belong to the shadow $S_o (x, R)$, then any $(Q,C)$-quasigeodesic segment between $y$ and $z$ does not intersect $S_o (x, R-D)$.
		\end{enumerate}
	\end{prop}
	
	\begin{proof}
		This result is known for geodesic hyperbolic spaces, \cite[Corollary 3.21]{maher_tiozzo18}. In the general case, due to Theorem \ref{thm isom embed geod hyp space}, there is an isometric embedding $f : (X, d) \rightarrow (Y, d_Y)$. Check that $y \in S_o(x, R) $ if and only if $f(y)$ belongs to the shadow $S^Y_{f(o)}(f(x), R) $ in $Y$. Applying \cite[Corollary 3.21]{maher_tiozzo18} on the shadow $S^Y_{f(o)}(f(x), R) $, we obtain the result. 
	\end{proof}
	
	We close this section with the notion of coarse metric. 
	\begin{Def}
		A \emph{coarse metric} on a set $X$ is a function $d: X \times X \to \bbR_+$ such that there exists $C>0$ for which the following conditions are verified for all $x,y,z \in X$: 
		\begin{enumerate}
			\item $d(x, x) \leq C$; 
			\item $d(x, y) \leq d(y, x) + C$; 
			\item $d(x, z ) \leq d(x, y) + d(y, z) + C$. 
		\end{enumerate}
		In this case, we say that $d$ is a \emph{$C$-coarse metric}.
	\end{Def}
	
	\subsection{Horofunction boundary}

	In this section, we show that the results of \cite{bader_caprace_furman_sisto22}, although written for geodesic hyperbolic spaces, can be adapted to the case of almost geodesic spaces. 
	
	Let $(X, d)$ be a separable metric space, and $o \in X$. Consider the horofunction map 
	 \begin{align} 
	 		\beta^o \colon  &X \rightarrow \bbR^X \nonumber \\
	 		&x \mapsto \beta_x^o : y \mapsto d(x,y) - d(x,o). \nonumber
	 \end{align}
	Let $\bbR^X$ be the set of all functions on $X$ with values in $\bbR$, and endow it with the topology of uniform convergence on compact sets, and consider the quotient topological vector subspace $\bbR^X / \bbR \cdot\mathds{1}$, with the quotient topology. We denote by $\Xhc$ the \emph{horoclosure} of $X$,i.e. the closure of $\beta^o (X)$ in $\bbR^X / \bbR\mathds{1}$, and by $\Xh$ the preimage of $\Xhc$ by the natural quotient $\bbR^X \rightarrow \bbR^X / \bbR\mathds{1}$. We will call \emph{horofunctions} the elements of $\Xhc$  (and by abuse of notations, of $\Xh$ too). We call $\partial^h X = \Xhc \setminus \beta^o (X)$ the \emph{horoboundary} of $X$. Since $X$ is assumed separable, $\Xh$ is metrizable and compact \cite[Lemma~3.2]{bader_caprace_furman_sisto22}. 
	
	We denote by $\Xhbc := \{h \in \Xhc \; | \; \inf (h) > -\infty\}$ and $\Xhuc := \{h \in \Xhc \; | \; \inf (h) = -\infty\}$. This decomposition is constant on the fibres in $\Xh$ so it gives a decomposition $\Xh = \Xhb \cup \Xhu$. Moreover, $\beta^o (X) \subseteq \Xhbc$, hence $\Xhbc$ is dense in $\Xhc$. We endow the spaces $\Xh$ and $\Xhc$ with the $\sigma$-algebra of Borel sets. The decompositions $\Xh = \Xhu \cup \Xhb$  and $\Xhc = \Xhbc \cup \Xhuc$ are measurable and $\iso(X) $-equivariant \cite[Lemma 3.3]{bader_caprace_furman_sisto22}. 
	
	For the rest of this section $(X, d)$ is a separable $\delta$-hyperbolic, $\alpha$-almost geodesic space. 
	The next Lemma is folklore, but most of the proofs available concern geodesic spaces so we provide one for self-containment. 
	
	\begin{lem}
		Let $(x_n)$ be a sequence of points in $X$ which converge in $\Xhc$ to the horofunction $\hat{\beta}$. Then $(x_n)$ converges to a point $\xi $ of the Gromov boundary $\bdg X$ if and only if $\hat{\beta} \in \Xhuc$. In this case, if $(x'_n)$ is another sequence such that $\lim x'_n = \hat{\beta}$, then $(x_n| x'_n )_o \rightarrow \infty$. 
	\end{lem}
	
	\begin{proof}
		Let $(x_n)$, $\hat{\beta}$ as in the Lemma, and denote by $b$ the lift of $\hat{\beta}$ in $\Xh_o$. We must show that $\beta \in \Xhu$. Assume that $(x_n)$ converges at infinity, that is $\underset{{n,m}}{\lim} (x_n|x_m)_o = \infty$. Fix $A > 0 $ and $N $ such that for $n, m > N$, $(x_n| x_m )_o \geq A$. For $m > N$, $(x_m|x_m)_o \geq A$ hence $d(x_m, o) \geq A$. By $\alpha$-almost geodesicity, take $z$ such that 
		\begin{align}
			|d(o,z)- A| \leq \alpha \text{ and } |d(z,x_m ) - (d(o,x_m) - A)|\leq \alpha. \nonumber
		\end{align}
		By hyperbolicity, for all $n \geq N$, 
		\begin{align}
			(x_n| z)_o \geq \min \{ (x_n|x_m)_o , (x_m| z)_o\} - \delta. \nonumber
		\end{align}
		But by definition, $d(x_n, z) - d(x_n, o)= d(o, z) - 2 (z|x_n)_o $, hence
		\begin{align}
			\beta(z) &= \lim_n (d(x_n, z) - d(x_n, o) ) \nonumber \\
			&\leq  d(o, z) + 2 \delta - 2(x_m| z)_o \nonumber \\
			& \leq  2 \delta + \alpha - A \nonumber. 
		\end{align}
		As $A$ was arbitrary, $\beta \in \Xhu_o$.
		
		The converse follows from \cite[Lemma 3.8]{maher_tiozzo18}: for any $\beta' \in \Xh_o$, and any $x, y \in X$, 
		\[\min\{- \beta'(x), -\beta'(y) \} \leq (x|y)_o + O(\delta). \]
		Using this inequality, one deduces the last part of the proof. 
	\end{proof}
	
	We can then define the map $\pi : \Xhuc \rightarrow \bdg X$ so that if $x_n \rightarrow \beta \in \Xhuc$,  $\pi (\lim_n x_n) = \xi  \in \bdg X$, where $\xi$ is the equivalence class of $(x_n)$. By \cite[Lemma 3.6]{bader_caprace_furman_sisto22}, this map is continuous and $\iso(X)$-equivariant. 
	
	Denote by $\bdd(X)$ the set of closed non-empty bounded subsets of $X$, and endow it with the Hausdorff metric and its Borel $\sigma$-algebra. 
	We consider the infimum function on $\Xhb$, defined by:
	\begin{align}
		\inf : h \in \Xhb \mapsto \inf\{h(x) \; | \; x \in X\} \nonumber. 
	\end{align}
	
	It is clearly $\iso(X)$-equivariant, and measurable. Indeed, fix a dense countable subset $X_0 \subseteq X$. By continuity of the functions in $\Xh$, $\inf(h) = \inf\{h(x) \; | \; x \in X_0\}$. 
	
	A key lemma in the proof of \cite[Theorem 3.1]{bader_caprace_furman_sisto22} is the following. 
	\begin{lem}[{\cite[Lemma~3.4]{bader_caprace_furman_sisto22}}]\label{lem coarse metric hyp}
		There is a Borel coarsely separable pseudo-metric on $\Xh$ on which $\iso(X)$ acts by isometries. Moreover, if for a discrete subgroup $G < \iso(X)$ the $G$-action on $X$ is unbounded, then so is the $G$-action on $\Xh$.
	\end{lem}
	The only part in the proof given in \cite{bader_caprace_furman_sisto22} of this lemma that uses the assumption that $X$ is geodesic is the following fact, which we reprove in our setting. 
	
	\begin{lem}
		Let $(X,d)$ be a $\delta$-hyperbolic, $\alpha$-almost geodesic separable metric space. For every $h \in \Xhb$, the set 
		\begin{align}
			\tilde{I} (h) :=\{x \in X \; | \; h(x ) < \inf(h) +1 \}\nonumber
		\end{align}
		is bounded. This gives a measurable, $\iso(X)$-equivariant map $\tilde{I} : \Xhb \rightarrow \bdd(X)$, and this map factors measurably and $\iso(X)$-equivariantly via $\Xhb \rightarrow \Xhbc$. 
	\end{lem}
	
	\begin{proof}
		Let $h \in \Xhb$, and assume without loss of generality that $\inf(h) = 0 $. Consider $x, x' \in \tilde{I}(h) $ and define $A :=d(x,x')$. As $X$ is $\alpha$-almost geodesic, take $z \in X$ such that 
		\begin{align}
			|d(x,z)- A/2| \leq \alpha \text{ and } |d(x', z ) - A/2 |\leq \alpha \nonumber.
		\end{align}
		Let $\hat{h}$ be the image of $h $ in $\Xhc$. By definition, there is a sequence $\{y_n\}$ in $X$ such that $\hat{h}$ is the image of the pointwise limit of elements of the form
		\begin{align}
			\beta_{y_n} : x \in X \mapsto d(x,y_n) \nonumber
		\end{align}
		in the quotient space $\bbR^X / \bbR\mathds{1}$. In particular, there exists $y \in X$ such that 
		
		\[  \left\{ \begin{array}{ll}
			|d(y,z) - h(z) - (d(y, x) - h(x))| < 1, \text{      and } \\
			|d(y,z) - h(z) - (d(y, x') - h(x'))| < 1.\end{array} \right. \] 
		
		Rewriting the previous system, this yields 
		
		\[  \left\{ \begin{array}{ll}
			h(z) < 1 + d(y , z ) + h(x) - d(y, x) \\
			h(z) < 1 + d(y , z ) + h(x') - d(y, x').\end{array} \right. \] 
		Since $h(x), h(x') < 1$, we obtain  
		\[  \left\{ \begin{array}{ll}
			h(z) < 2 + d(y , z ) - d(y, x) \\
			h(z) < 2 + d(y , z ) - d(y, x').\end{array} \right. \] 
		
		Fix such a $y \in X$. By Gromov hyperbolicity:
		\begin{align}
			(x,x')_y \geq \min \{(x,z)_y, (x', z)_y\} - \delta \nonumber. 
		\end{align}
		Rewriting this, we have two possibilities. If $(x,z)_y \geq (x', z)_y$, 
		\begin{align}
			&d(x,y) + d(x', y) - d(x, x') \geq d(x', y) + d(z, y) - d(x', z) -2 \delta \nonumber \\
			\Leftrightarrow & d(z, y ) \leq d(x,y) + d(x', z) - d(x,x') + 2 \delta. \label{eq geq}
		\end{align}
		If $(x,z)_y \leq (x', z)_y$, the same computation gives
		\begin{align}
			d(z, y ) \leq d(x',y) + d(x, z) - d(x,x') + 2 \delta. \label{eq leq}
		\end{align}

		Now if $(x,z)_y \geq (x', z)_y$, we have in particular by equation \eqref{eq geq}:
		\begin{align}
			h(z) &< 2 + d(y , z ) - d(y, x) \nonumber\\
			& \leq 2 + d(x,y) + d(x', z) - d(x,x') + 2 \delta - d(x,y) \nonumber \\
			&\leq 2 + 2 \delta + d(x', z ) - d(x,x') \nonumber. 
		\end{align}
		If $(x,z)_y \leq (x', z)_y$, equation \eqref{eq leq} implies 
		\begin{align}
			h(z) < 2 + 2 \delta + d(x, z) - d(x,x') \nonumber
		\end{align}
		
		But by definition of $z$, 
		\begin{align}
			d(x, z) - d(x,x') \leq \alpha - A/2  \text{ and } d(x', z) - d(x,x') \leq  \alpha - A/2 \nonumber
		\end{align}
		Hence, taking $A > 2\alpha + 4 + 4 \delta $ gives $h(z) < 0 $, a contradiction with the fact that $\inf(h) = 0$. Then $A$ is (uniformly) bounded. 
		
		Starting from there, the arguments of \cite[Lemma~3.4]{bader_caprace_furman_sisto22} do not use the assumption that $X$ is geodesic and we can proceed in the exact same way.  
	\end{proof}

	We close this section with some terminology. 
	Recall that the \emph{support} of a Borel measure $m$ on a topological space $Y$ is the smallest closed set $C$ such that $m(Y \setminus C)= 0$. In other words $y \in \supp(m)$ if and only if for all $U $ open containing $y $, $m(U) >0$.
	
	\begin{Def}
		We say that the action by isometries of a group $G$ on a hyperbolic space $Y$ (not assumed to be proper) is \emph{non-elementary} if there are two loxodromic isometries with disjoint fixed points on the Gromov boundary. In that case, a probability measure $\mu$ on $G$ is said to be \emph{non-elementary} if its support generates a group acting non-elementarily on $Y$. 
	\end{Def}
	For $G < \iso(X)$ a group of isometries, we denote by 
	\[\Lambda(G):= \{ \eta \in \bdg X \mid g_n x \to \eta \text{ for some } (g_n) \in G^\bbN\}\]
	the \emph{limit set} of $G$. 
	By the classification theorem of isometric actions on hyperbolic spaces \cite[Theorem~6.2.3]{das_simons_urbanski17}, a (semi)group of isometries of a hyperbolic space is non-elementary if and only if $\Lambda(G) $ is an uncountable set and $G$ does not have a global fixed point, which is the case if and only if there is no bounded orbit in $X$ nor finite orbit on $\bdg X$.  
\subsection{On the separability assumption}\label{section separable}
	Throughout the paper, we assume that the roughly geodesic hyperbolic metric space $(X, d)$ is separable. Since we only consider discrete countable groups acting on $(X, d)$, this is harmless.
	Indeed, we can pass to an action on a separable hyperbolic metric space $Y$ for which there is an equivariant quasi-isometric embedding $Y \hookrightarrow X$, by considering countably many orbits and attaching rough geodesics equivariantly. The procedure is explained in \cite[Remark~4]{gruber_sisto_tessera17} for geodesic spaces, but works in the exact same way for roughly geodesic almost geodesic hyperbolic spaces.

\section{Hyperbolic models for $\cat$(0) spaces}\label{section hyperbolic models}

	The goal of this section is to briefly present some ideas of \cite{petyt_spriano_zalloum24}, in which the authors build a way of attaching a family of hyperbolic metric spaces $X_L = (X, d_L)_L$ to a proper $\cat$(0) space. What is interesting about these spaces is that they convey much of the geometry of the original space, especially at infinity, and they behave very well under isometric actions. Moreover, contracting isometries act on these hyperbolic models as loxodromic isometries. Throughout this section, $(X,d)$ is a $\cat(0)$ space. 
	
	\subsection{Curtains and chains}
	An important feature of $\cat(0)$ spaces is the existence of closest-point projections on complete convex subsets \cite[Lemma II.2.4]{bridson_haefliger99}: given a complete convex subset $C\subseteq X$ in a $\cat$(0) space, there exists a map $\pi_C : X \rightarrow C$ such that for all $x \in X$, $\pi_C(x)$ minimizes the distance $d(x,C)$. For $I$ a closed interval, if $\gamma : I \rightarrow X$ is a geodesic, we will write $\pi_\gamma$ for the closest-point projection onto the image of $\gamma $ in $X$. 
	
	\begin{Def}
		Let $X$ be a $\cat(0) $ space, and let $\gamma: I \rightarrow X$ be a geodesic. Let $t \in I$ be such that $[t - \frac{1}{2} , t + \frac{1}{2}]$ belongs to $I$. Then the \emph{curtain} dual to $\gamma$ at $t$ is 
		\begin{align}
			h = h_{\gamma, t} = \pi^{-1}_\gamma (\gamma ([t - \frac{1}{2} , t + \frac{1}{2}])). \nonumber
		\end{align}
		The sets $h^{-} =  \pi^{-1}_\gamma (\gamma ((-\infty, t- \frac{1}{2})\cap I))$ and $h^{+} =  \pi^{-1}_\gamma (\gamma ((t +\frac{1}{2}, + \infty )\cap I))$ are the \emph{halfspaces} determined by $h$. Note that $\{h^{-} , h , h^{+}\}$ is a partition of $X$. If $A \subseteq h^{-}$ and $B \subseteq h^{+}$ are subsets of $X$, we say that $h$ \emph{separates} $A$ from $B$. 	
	\end{Def}

	A family of curtains $\{h_i\}$ is said to be a chain if $h_i $ separates $h_{i-1}$ from $h_{i+1}$ for every $i $. Chains can be used in order to define a metric on $X$ by the following: for $x \neq y \in X$, 
	\begin{align}
		d_\infty (x, y) = 1 + \max \{\, |c| \, : \, c \text{ is a chain separating }x \text{ from } y \}. \nonumber
	\end{align}
	
	If $h$ be a curtain, then $d(h^{-}, h^{+}) = 1$, hence for any $x,y \in X$, $d_\infty (x, y ) \leq \lceil d(x,y ) \rceil$. Conversely, $d$ and $d_\infty$ differ by at most 1 \cite[Lemma 2.10]{petyt_spriano_zalloum24}.
	
	We say that a chain $c$ of curtains \emph{meets} a curtain $h$ if every single curtain $h_i \in c$ intersects $h$. 
	
	\begin{Def}[$L$-separation]
		Let $L \in \mathbb{N}^\ast$, we say that disjoint curtains are $L$\emph{-separated} if every chain meeting both has cardinality at most $L$. A chain of pairwise $L$-separated curtains is called an $L$-\emph{chain}. 
	\end{Def}
	
\begin{figure}
		\centering
		\begin{center}
			\begin{tikzpicture}[scale=1.2]
				\draw (0,0) -- (4,0)  ;
				\draw (2.8,2) -- (2.8,-1)  ;
				\draw (1.2,2) -- (1.2,-1)  ;
				\draw (2,2) -- (2,-1)  ;
				\draw (2,0.5) node[right]{$\leq 2L+1$} ;
				\draw (-1, 2.5) to[bend right = 80] (5, 3);
				\draw (0,0) node[below left]{$x_1$} ;
				\draw (2, -1) node[below]{$h_2$} ;
				\draw (1.2, -1) node[below]{$h_1$} ;
				\draw (2.8, -1) node[below]{$h_3$} ;
				\draw (-1, 2.5) node[above left]{$x_2$} ;
				\draw (4,0) node[below right]{$y_1$} ;	
				\draw (5,3) node[above]{$y_2$};
				\filldraw[black] (2,0) circle(1pt);
				\filldraw[black] (2,1) circle(1pt);
			\end{tikzpicture}
		\end{center}
		\caption{Illustration of Lemma \ref{lem bottleneck}.}\label{bottleneck figure}
	\end{figure}
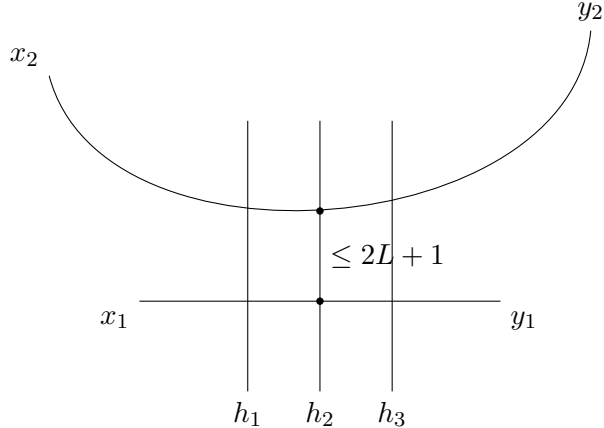

	The following geometric lemma means that $L$-separation induces good Morse properties, see Figure \ref{bottleneck figure}. 
	
	\begin{lem}[{Bottleneck Lemma, \cite[Lemma 2.14]{petyt_spriano_zalloum24}}]\label{lem bottleneck}
		Suppose that $A$, $B$ are two sets which are separated by an L-chain $\{h_1, h_2, h_3\}$ all of whose elements are dual to a geodesic $\gamma= [x_1, y_1]$ with $x_1 \in A$ and $y_1 \in B$. Then for any $x_2 \in A$, $y_2 \in B$, if $p \in h_2 \cap [x_2, y_2]$, then $d(p, \pi_\gamma(p))\leq 2L + 1$. 
	\end{lem}
	
	If there is a $L$-chain separating two points $x $ and $y$, we can find a (smaller) $L$-chain of curtains separating these, which is dual to the geodesic $[x,y]$ and whose size can be controlled. 
	
	\begin{lem}[{\cite[Lemma 2.21]{petyt_spriano_zalloum24}}]\label{lem dual chain}
		Let $L, n \in \mathbb{N}$, and let $\{h_1, \dots, h_{(4L +10)n}\}$ be an $L$-chain separating $A$, $B \subseteq X$. Take $x \in A$, $y \in B$. Then $A$ and $B$ are separated by an $L$-chain of size $\geq n+1 $ dual to $[x,y]$. 
	\end{lem}
	
\subsection{Hyperbolic models and isometric equivariance}
 	
	Given distinct points $x\neq y \in X$, we define 
	\begin{align}
		d_L(x,y) = 1 + \max\{|c| \, : \, c \text{ is an } L\text{-chain separating }x \text{ from } y\}. \nonumber
	\end{align}
	
	For every $L$, $d_L$ gives a metric on $X$ \cite[Lemma~2.17]{petyt_spriano_zalloum24}. We will denote by $X_L = (X, d_L)$ the resulting metric space. The metric spaces $(X, d_L)$ are quasi-geodesic hyperbolic space with hyperbolicity constants depending only on $L$, and $\iso(X)$ acts by isometries on $(X, d_L)$ \cite[Theorem~3.9]{petyt_spriano_zalloum24}. 
	Moreover, a semisimple isometry of $(X, d) $ is contracting if and only if there exists $L \in \mathbb{N}$ such that the corresponding isometry on $(X, d_L)$ is loxodromic \cite[Theorem~4.9]{petyt_spriano_zalloum24}.

\subsection{Boundaries}

	\begin{Def}
		We say that a geodesic ray $\gamma : [0, \infty) \rightarrow X$ \emph{crosses} a curtain $h $ if there exists $t_0 \in [0, \infty)$ such that $h$ separates $\gamma(0) $ from $\gamma ([t_0, \infty))$. Alternatively, we may say that $h$ separates $\gamma(0)$ from $\gamma(\infty)$. Similarly, we say that a geodesic line $\gamma: \mathbb{R} \rightarrow X$ crosses a curtain $h$ if there exist $t_1, t_2 \in \mathbb{R}$ such that $h$ separates $\gamma ((-\infty, t_1])$ from $\gamma ([t_2, \infty))$. We say that $\gamma $ \emph{crosses a chain} $c = \{h_i\}$ if it crosses each individual curtain $h_i$. 
	\end{Def}
	
	As a consequence of Lemma \ref{lem bottleneck} and Lemma \ref{lem dual chain}, if two geodesic rays with the same starting point cross an infinite $L$-chain $c$, then they are asymptotic, and hence equal. 

	Given $o \in X$, we define $\mathcal{B}_L $ as the subspace of $\bd X$ consisting of all geodesic rays $\gamma : [0, \infty) \rightarrow X$ starting from $o$ and such that there exists an infinite $L$-chain crossed by $\gamma$.
	A topology on $B_L$ can be given by the following basis. For a geodesic ray $\gamma$ emanating from $o $ such that $\gamma(\infty) = \xi \in \bd X$, and $h $ a curtain dual to $\gamma$, define 
	\begin{align}
		U_h(\xi) := \{ \eta \in \bd X \; | \; \gamma_o^\eta \text{ crosses $h$}\}, \nonumber
	\end{align}
	where $\gamma_o^\eta$ is the unique geodesic ray based at $o$, in the class of $\eta$. 
	The usual visual topology on $\bd X$ is denoted by $\mathcal{T}_{\text{cone}}$, see for instance \cite[Chapter II.8]{bridson_haefliger99}.
	We say that $U \subseteq \bd X$ is open for the \emph{curtain topology} if for every $\xi \in U$, there exists a curtain $h $ dual to $\gamma_o^\xi $ such that $U_h(\xi) \subseteq U$. The identity map 
	\begin{align}
		(\bd X, \mathcal{T}_{\text{cone}}) \rightarrow (\bd X,\mathcal{T}_{\text{curtain}}) \nonumber
	\end{align} 
	is continuous, and the topologies are the same in restriction to any $B_L$ \cite[Theorem 8.8]{petyt_spriano_zalloum24}.

	In the case of the contact graph associated to a $\cat$(0) cube complex $X$, we had the existence of an $\iso(X)$-equivariant embedding of the boundary of the contact graph into the Roller boundary $\partial_{\mathcal{R}}X$ \cite[Theorem 3.7]{fernos_lecureux_matheus21}.  The following analogous result was proven in the context of proper $\cat$(0) spaces by Petit, Spriano and Zalloum \cite[Theorem 8.1]{petyt_spriano_zalloum24}. 
	
	\begin{thm}\label{thm homeomorphism bords XL}
		Let $X$ be a complete $\cat$(0) space. Then for any $L\geq 0 $, the identity map $\iota_L : (X,d) \rightarrow (X,d_L) $ induces an $\iso(X)$-equivariant homeomorphism of the boundary $\partial_L : B_L \rightarrow \bdg X_L$. 
	\end{thm}

	We begin with a lemma. 
	\begin{lem}[{\cite[Corollary 8.11]{petyt_spriano_zalloum24}}]\label{lem dual chain qgeod}
		Let $P : [0, \infty) \rightarrow X$ be a $(Q,C)$-quasigeodesic ray of $X_L $. Then there is a sequence $(x_i)_{i \in \mathbb{N}} \subseteq P $ and a $L$-chain $\{c_i\}_{i \in \mathbb{N}}$ such that for every $n \in \mathbb{N}$, $\{c_1, \dots, c_n\} $ separates $o $ from $x_n $. 
	\end{lem}
	
	The following proposition will be essential. 
	
	\begin{prop}\label{prop inverse map}
		Let $(x_n) $ be a sequence in $X$ such that $(x_n)$ converges to a point $\xi_L$ in the Gromov boundary $\bdg X_L$ (for the hyperbolic topology). Then $(x_n)$ converges to a point $\xi \in \bd X$ in the visual boundary of $X$. 
	\end{prop}
	
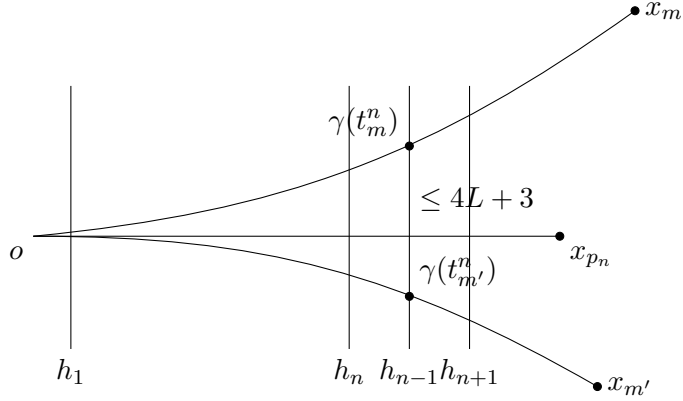
\begin{figure}
		\centering
		\begin{center}
			\begin{tikzpicture}[scale=1]
				\draw (-3,0) -- (4,0)  ;
				\draw (2.8,2) -- (2.8,-1.5)  ;
				\draw (-2.5,2) -- (-2.5,-1.5)  ;
				\draw (1.2,2) -- (1.2,-1.5)  ;
				\draw (2,2) -- (2,-1.5)  ;
				\draw (2,0.5) node[right]{$\leq 4L+3$} ;
				\draw (-3, 0) to[bend right = 15] (5, 3);
				\draw (-3, 0) to[bend left = 15] (4.5, -2 );
				\draw (-3,0) node[below left]{$o$} ;
				\draw (2, -1.5) node[below]{$h_{n-1}$} ;
				\draw (1.2, -1.5) node[below]{$h_n$} ;
				\draw (-2.5, -1.5) node[below]{$h_1$} ;
				\draw (2.8, -1.5) node[below]{$h_{n+1}$} ;
				\draw (4,0) node[below right]{$x_{p_n}$} ;	
				\draw (2,-0.8) node[above right]{$\gamma(t^n_{m'})$} ;
				\draw (2,1.2) node[above left]{$\gamma(t^n_{m})$} ;
				\draw (5,3) node[right]{$x_{m}$};
				\draw (4.5, -2 ) node[right]{$x_{m'}$};
				\filldraw[black] (2,-0.8) circle(1.5pt);
				\filldraw[black] (2,1.2) circle(1.5pt);
				\filldraw[black] (4.5,-2) circle(1.5pt);
				\filldraw[black] (4,0) circle(1.5pt);
				\filldraw[black] (5,3) circle(1.5pt);
			\end{tikzpicture}
		\end{center}
		\caption{Proof of Proposition \ref{prop inverse map}.}\label{figure prop inverse map}
	\end{figure}

	\begin{proof}[Proof of Proposition \ref{prop inverse map}]
		As $X_L$ is a $\alpha $-almost geodesic space, Proposition \ref{prop exist qgeod hyp} implies that $X_L$ is a $(Q,C) $-quasigeodesic space. In particular there exists a $(Q,C)$-quasigeodesic $P : [0 , \infty) \rightarrow X_L $ emanating from $o$ and representing $\xi_L$. By Lemma \ref{lem dual chain qgeod}, there exists a sequence $(y_i) $ in $P$ and an infinite $L$-chain $\{c_i\}$ such that for every $n \in \mathbb{N}$, $\{c_1, \dots, c_n\} $ separate $o $ from $ y_n$. Since $x_n \rightarrow \xi_L$, for every $n \in \mathbb{N}$, there exists $p_n \geq 0 $ such that for all $m \geq p_n $, $\{c_0, \dots, c_{(4L + 10)n} \} $ is a $L$-chain separating $o $ from $x_m$. 
		
		Now, define the geodesic segments $\gamma_n = [o, x_n]$ (for the $\cat$(0) distance) parametrized by arc-length. In other words for every $t \in [0, \infty)$, $\gamma_n(t) $ is the point at distance $t$ from $o$ on the geodesic segment $[0, x_n]$ if $t \leq d(o, x_n)$ and $\gamma_n (t) = x_n$ if $t \geq d(o, x_n)$. We are going to prove that for every $R> 0$, $\gamma_n(R)$ is a Cauchy sequence. 
		
		Let $R>0$ be fixed. By Lemma \ref{lem dual chain}, there exists a $L$-chain of size $n+1$, dual to the geodesic segment $[o, x_{p_n}]$, and separating $c_1 $ from $c_{(4L+10)n}$. Denote this chain by $\{h_1, \dots, h_{n+1}\}$. Since curtains are thick, the number of curtains in $\{c'_i\}$ that intersect the metric ball $B(o, R) $ is less than $R+1$. In particular, for $n $ large enough and for every $m \geq p_n$, the number of curtains in $\{h_1, \dots, h_{n+1}\}$ separating $\gamma_m (R)$ and $x_m $ is greater than $n -R$. Let $n \geq R + 3$, so that $\{h_{n-R+1}, \dots, h_{n+1}\}$ is a $L$-chain separating $B(o, R)$ from $x_m $ for every $m \geq p_n$. 
		
		Let $m \geq p_n$, and take $t^n_m \in [0, d(o, x_m)]$ such that $\gamma_m(t^n_m) \in h_n$, which is possible since $\gamma_m$ crosses $h_n$, see Figure \ref{figure prop inverse map}. By the Bottleneck Lemma \ref{lem bottleneck}, for all $m , m'\geq p_n$, $d(\gamma_m(t^n_m), \gamma_{m'}(t^n_{m'})) \leq 4L+3$. Without loss of generality, we can assume that $t^n_m \leq t^n_{m'} $. Note that because curtains are thick, and since the projections do not increase distances (Proposition \cite[Lemma~II.2.4]{bridson_haefliger99}), $d(o, \gamma_m(t_m^n)) \geq d(o, h_n)$, hence $t^n_m \geq n$. By convexity of the $\cat$(0) metric \cite[Proposition II.2.2]{bridson_haefliger99}, we have 
		\begin{align}
			d(\gamma_m(t^n_m), \gamma_{m'}(t^n_{m})) \leq 4L+3. \nonumber
		\end{align}
		
		The function $t \rightarrow d(\gamma_m(t), \gamma_{m'}(t))$ is convex, hence 
		\begin{align}
			d(\gamma_m(R), \gamma_{m'}(R)) &\leq \frac{R}{t^n_m} d(\gamma_m(t^n_m), \gamma_{m'}(t^n_{m}))\nonumber \\
			& \leq  \frac{R}{t^n_m} (4L+3) \nonumber \\
			& \leq  \frac{R}{n} (4L+3) \nonumber. 
		\end{align}
		
		As a consequence, for every $\varepsilon >0$, there exists $p_\varepsilon \in \mathbb{N}$ large enough such that for all $m, m' \geq p_n$, $d(\gamma_m(R), \gamma_{m'}(R)) \leq \varepsilon$. Therefore, for all $R >0$, $\gamma_m(R)$ is a Cauchy sequence, and because $X$ is complete, the sequence $(\gamma_n)_n$ converges to a geodesic ray. 
	\end{proof}
	
	The rest of the proof of Theorem \ref{thm homeomorphism bords XL} is now very similar to what was done in \cite[Proposition 8.9]{petyt_spriano_zalloum24}. We include it for clarity. 	
	
	\begin{proof}[Proof of Theorem \ref{thm homeomorphism bords XL}]
		The map $\iota_L : (X,d) \rightarrow (X,d_L)$ is $(1,1)$-coarsely Lipschitz, hence by \cite[Lemma 6.18]{incerti-medici_zalloum23}, the existence and continuity of $\partial_L$ follows from the fact that $X_L$ is $\alpha$-almost geodesic.

		If $\gamma$ crosses an infinite $L$-chain, then there is an infinite $L$-chain $\{h_i\}$ dual to $\gamma$. Orient this chain so that $o \in h^{-1}_i $ for all $i$. Let $\gamma '$ be a distinct geodesic ray, emanating from $o $. Then by Lemma \ref{lem bottleneck} and there must be a $k $ such that $\gamma' \subseteq h_k^{-1}$. In particular, if $x_n \in h_n \cap \gamma$, then $d_L(\gamma', x_n) \geq n-k$, and $\gamma'$ does not lie in a finite $X_L$ neighbourhood of $\gamma$. This argument shows that the map $\partial_L$ is injective.

		Let $P : [0, \infty) \rightarrow X_L$ be a quasigeodesic ray (for the metric $d_L$) emanating from $o$, and let $(x_n)$ be an unbounded sequence on $P$. By Lemma \ref{lem dual chain qgeod}, there exists an infinite $L$-chain $\{c_i\}$ such that for every $n \in \mathbb{N}$, $\{c_1, \dots, c_n\} $ separate $o $ from $ x_n$. Proposition \ref{prop inverse map} shows that $[o, x_n]$ converges to a geodesic ray $\gamma$. By construction, $\gamma$ crosses infinitely many curtains $\{c_i\}$, hence belongs to $B_L$. Note that by Proposition \ref{cor stab qgeod hyp}, any $(1,C)$-quasigeodesic ray representing the same endpoint $\xi_L$ as $P$ lies in a uniform neighbourhood of $P$, with universal constants. Therefore, $\iota(\gamma)$ is an unparametrized rough geodesic, and lies in a uniform neighbourhood of any quasigeodesic ray between $\iota(o)$ and $\xi_L$. Then, the choice of $\gamma$ does not depend upon the choice of $P$, and the application $\partial_L$ is surjective.

		It remains to prove that the inverse map is continuous. This part is completely similar to what is presented in $(4)$ of the proof of \cite[Proposition 8.9]{petyt_spriano_zalloum24}, which does not use any properness assumption on $X$. 
	\end{proof}

	As a consequence, we obtained the following proposition.  
	
	\begin{prop}\label{prop non elem}
		Let $G$ be a group acting non-elementarily by isometries on a complete $\cat$(0) space $(X,d)$, and assume that $G$ possesses a pair of independent contracting elements for this action. Then there exists $L \in \mathbb{N}$ such that $G$ acts on the hyperbolic space $(X, d_L)$ non-elementarily by isometries. 
	\end{prop}
	
	\begin{proof}
		By \cite[Theorem~4.9]{petyt_spriano_zalloum24}, there exists $L \in \mathbb{N}$ such that $g$ and $h$ act on $(X, d_L)$ as loxodromic isometries. As $g $ and $h$ are independent, their fixed points form four distinct points of the visual boundary $\bd X$. Now seen in $X_L = (X, d_L)$, their fixed points sets must also form four distinct points of $\partial X_L$ because of the homeomorphism $\partial_L : \mathcal{B}_L \longrightarrow \partial X_L$. This means that the action $G \curvearrowright X_L$ is non-elementary. 
	\end{proof}

	\section{Random dynamical systems}\label{section rds}
	
	This section introduces the notations and first properties of random dynamical systems. 
	
	\subsection{RDS and greg's}
	
	A standard Borel space is a measurable space measurably isomorphic to a Polish space. A standard measured (resp. probability) space is a standard Borel space equipped with a Borel measure (resp. probability measure). Let $(\Omega, \mathcal{F}, \bbP)$ be a standard probability space. Let $\mathbb{T}= \mathbb{Z}$ (two-sided time) or $\mathbb{N}$ (one-sided time) and let $\{\theta (t)\}_{t \in \mathbb{T}}$ be a family of measure-preserving transformations on  $(\Omega, \calF, \bbP)$ satisfying 
	$$ \theta (t)\theta (s)= \theta (t+s) \text{ for all $s,t$ and } \theta(0)=\Id.$$
	The family $\{\theta (t)\}_{t \in \mathbb{T}}$ can be called a \emph{flow} on $\Omega$. Let $X$ be a Polish space, with Borel $\sigma$-algebra $\mathcal{B}$. A \emph{random dynamical system} on $(X, \mathcal{B})$ over $(\Omega, \mathbb{P}, \{\theta (t)\})$ is a measurable map 
	\begin{align}
		\varphi \colon& \mathbb{T} \times \Omega \times X \to X \nonumber \\
		& (t, \omega, x) \mapsto \varphi(t, \omega)x \nonumber 
	\end{align}
	such that $\varphi(0, \omega)x = x $ for all $(\omega, x)$ and such that for all $s, t \in \mathbb{T}$ and for all $\omega $ in a $\{\theta (t)\}$-invariant set of full $\mathbb{P}$-measure, we have the cocycle relation 
	\begin{align}\label{eq cocy rds}
		\varphi(t+s, \omega) = \varphi(t, \theta(s) \omega)\circ\varphi(s, \omega). 
	\end{align}
	In fact, we could have taken the same definition for a \emph{continuous flow}, i.e. when $\mathbb{T}= \bbR$ or $\bbR_+$, but we will only deal with discrete time. 
	
	\begin{rem}
		Sometimes, RDS of the previous form are called \emph{perfect}, because there is a $\mathbb{P}$-conull set on which the cocycle relation \eqref{eq cocy rds} is verified for all $t,s \in \mathbb{T}$. When the set of $\omega $ for which the cocycle relation holds depends on the time $s$, the RDS is called \emph{crude}. Crude cocycles can be made perfect \cite[Chapter 1.3]{arnold98}, so all the RDS we shall consider are assumed perfect. 
	\end{rem}

	If the flow $\{\theta (t)\}$ is two-sided ($\mathbb{T}= \mathbb{Z}$), we immediately have that 
	$$ \varphi(t, \omega)^{-1} = \varphi(-t, \theta(t)\omega).$$
	Given $\varphi$ a RDS on $(X, \mathcal{B})$ over $\{\Omega, \mathbb{P}, \{\theta (t)\}\}$, we let $\check{\varphi}$ and $\chi$ be the measurable maps defined by 
	$$\check{\varphi}(t, \omega) = \varphi(t, \omega)^{-1} = \varphi(-t, \theta(t)\omega)$$
	and 
	$$\chi(t, \omega) = \varphi(-t, \omega)^{-1} = \varphi(t, \theta(-t)\omega)$$
	for all $ t \in \mathbb{T}$ and almost every $\omega $. 
	
	The maps  $\check{\varphi}$ and $\chi$ satisfy the ``backward'' cocycle relations: for all $s, t \in \mathbb{T}$ and almost every $\omega $, 
	$$ \check{\varphi}(t+s, \omega) = \check{\varphi}(s, \omega) \circ \check{\varphi}(t, \theta(s)\omega)$$ and 
	$$ \chi(t+s, \omega) = \chi(s, \omega) \circ \chi(t, \theta(-s)\omega).$$
	
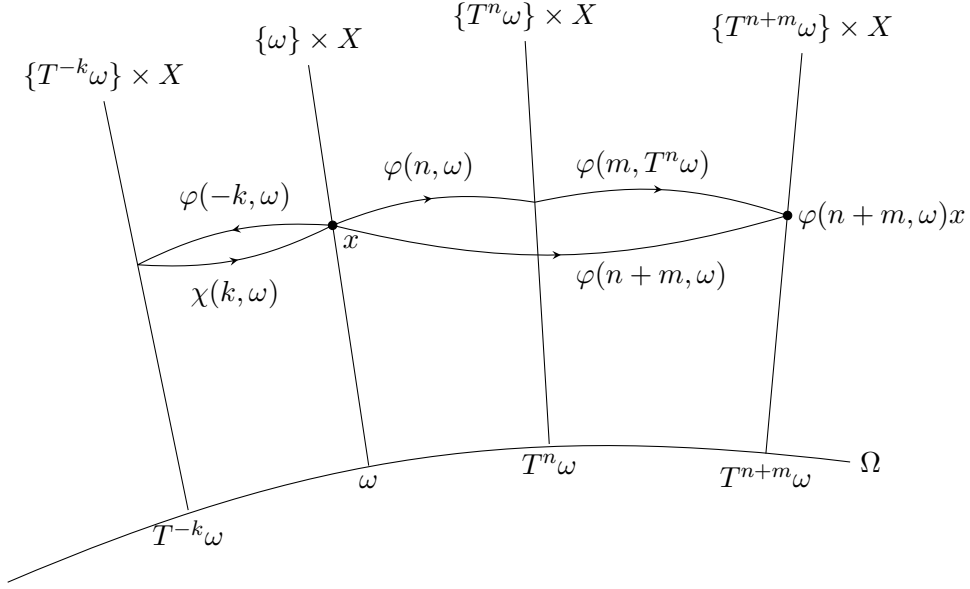
\begin{figure}[h]
		\centering
		\begin{center}
			\begin{tikzpicture}[scale=1.6]
				\draw (-3, 0) to[bend left = 15] (4, 1);
				\draw (-1.5,0.6) --  (-2.2,4);
				\draw (0,0.97) --  (-0.5,4.3);
				\draw (1.5,1.15) --  (1.3,4.5);
				\draw (3.3,1.07) --  (3.6,4.4);
				\draw (-1.5,0.6) node[below]{$T^{-k} \omega$} ;
				\draw (0,0.97) node[below]{$ \omega$} ;
				\draw (1.5,1.15) node[below]{$T^{n} \omega$} ;
				\draw (3.3,1.07) node[below]{$T^{n+m} \omega$} ;
				\draw (-1.11,3) node[above]{$\varphi(-k, \omega)$} ;
				\draw (0.5,3.25) node[above]{$\varphi(n, \omega)$} ;
				\draw (2.28, 3.27) node[above]{$\varphi(m, T^n\omega)$} ;
				\draw (2.35, 2.75) node[below]{$\varphi(n+m, \omega)$} ;
				\draw (-1.11, 2.6) node[below]{$\chi(k, \omega)$} ;
				\draw (-0.5,4.3) node[above]{$\{\omega\}  \times X$} ;
				\draw (-2.2,4) node[above]{$\{T^{-k}\omega\}  \times X$} ;
				\draw (1.3,4.5) node[above]{$\{T^{n}\omega\}  \times X$} ;
				\draw (3.6,4.4) node[above]{$\{T^{n+m}\omega\}  \times X$} ;
				\draw (-0.3, 2.968) node[below right]{$x$} ;
				\draw (3.48,3.05) node[right]{$\varphi(n+m, \omega)x$} ;
				\draw (4, 1) node[right]{$\Omega$} ;
				\filldraw[black] (-0.3, 2.968) circle(1pt);
				\filldraw[black] (3.48,3.05) circle(1pt);
				\path [draw=black,postaction={on each segment={mid arrow=black}}]
				(-1.92, 2.64) to[bend right = 15]  (-0.3, 2.968)
				(-0.3, 2.968) to[bend right = 15]  (-1.92, 2.64)
				(-0.3, 2.968) to[bend left = 15] (1.38,3.16)
				(1.38,3.16) to[bend left = 15] (3.48,3.05)
				(-0.3, 2.968) to[bend right = 15] (3.48,3.05)
				;
				
				%\draw (2, -1.5) node[below]{$h_{n-1}$} ;
				%\filldraw[black] (2,-0.8) circle(1.5pt);
			\end{tikzpicture}
		\end{center}
		\caption{A discrete RDS seen as a bundle $\Omega \times X$ over $\Omega$.}\label{figure rds}
	\end{figure}

	We say that $\check{\varphi}$ is a \emph{backward cocycle} with respect to the flow $\{\theta (t)\}$. Accordingly $\chi $ is a backward cocycle with respect to the flow $\{\theta^{-1} (t)\}$. It is useful to visualize these cocycles as if they are acting on the bundle $\Omega \times X$ over $\Omega$, see Figure \ref{figure rds}. 
	
	\begin{rem}
		We use the term ``backward'' as in \cite{arnold98}: it means that the cocycles $\check{\varphi}$ and $\chi$ behave in a somewhat acausal manner, especially when associated with group actions.
	\end{rem}
	
	We will consider discrete RDS of the following particular form. Let $(\Omega, \mathcal{F}, \mathbb{P})$ be a standard probability space, and $T : \Omega \to \Omega$ be a measure preserving and invertible transformation. Let $G$ be a locally compact second countable group, and let $f: \Omega \to G$ be a measurable map. 
	\begin{Def}
		A \emph{group random element generator} (greg) $(\Omega, \mathbb{P}, T, G, f)$ is the invertible measurable map defined by 
		\begin{align}
			\varphi \colon &\mathbb{Z} \times \Omega \to G \nonumber \\
			& (t, \omega) \mapsto \varphi(t, \omega), \nonumber 
		\end{align}
		with $\varphi(1, \omega)= f(\omega)$ and  for all $t,s \in \mathbb{Z}, \omega \in \Omega$,
		$$\varphi(t+s, \omega) = \varphi(t, T^s\omega)\varphi(s,\omega).$$
	\end{Def}
	
	Note that in this case, by the cocycle relation \eqref{eq cocy rds}, we have that 
	\begin{align}
		&\Id = \varphi(0, \omega)= \varphi(1 -1, \omega) = \varphi(1, T^{-1}\omega)\varphi(-1, \omega) \nonumber \\
		\Rightarrow \ & \varphi(-1, \omega) = \varphi(1, T^{-1}\omega)^{-1} = f(T^{-1}\omega)^{-1}. \nonumber
	\end{align}	
	
	As a consequence, the greg $(\Omega, \mathbb{P}, T, G, f)$ can be completely described by $\varphi(0, \omega) = \Id $ and for all $n \geq 1, \omega \in \Omega$,
	\begin{align}
		\varphi(n, \omega) &= f(T^{n-1}\omega)\dots f(\omega) \nonumber \\
		\varphi(-n, \omega) &= f(T^{-n}\omega)^{-1}f(T^{-n +1}\omega)^{-1}\dots f(T^{-1}\omega)^{-1} \nonumber.
	\end{align}

	The cocycle $(\varphi(n, \omega))_{n\in \mathbb{Z}}$ describes a stochastic walk where the increments are distributed according to $f_\ast \mathbb{P}$, but are not necessarily independent. If $G$ acts on the standard Borel space $(X, \calB)$ by Borel maps, the associated RDS over  $(\Omega, \mathbb{P}, T, G, f)$ is :
	\begin{align}
		\varphi \colon& \mathbb{Z} \times \Omega \times X \to X \nonumber \\
		& (t, \omega,x) \mapsto \varphi(t, \omega)x.\nonumber 
	\end{align}

\subsection{Disintegration of measures}%\label{section disintegration}
	Let $(Y, \mathcal{Y})$ be a standard Borel space and let $\mu$ be a probability measure on $(Y, \mathcal{Y})$. Let $p : (Y, \mathcal{Y}) \to (Z, \mathcal{Z})$ be a Borel map, and denote by 
	\[\mathcal{Y}' = p^{-1}\mathcal{Z} = \{ p^{-1} E \mid E \in \mathcal{Z}\}\] the sub-$\sigma$-algebra of $\mathcal{Y}$ given by this map. The disintegration Theorem of Rohlin \cite{rholin49} gives the following: there exists a $ \mu$-essentially unique $\mathcal{Y}'$-measurable map $y \mapsto \mu_y \in \prob(Y)$ such that, for $\mu$-almost every $y \in Y$, $\mu_y$ is supported by $p^{-1}(p(y))$ and
	$$ \mu = \int_Y \mu_y d\mu(y).$$ 
	For any non-negative measurable function $h : Y \to \bbR_+$, we have 
	$$ \mathbb{E}(h| \mathcal{Y}') (y) = \int_Y h(y') d\mu_y(y').$$
	
	The measurable map $y \mapsto \mu_y$ is called the \emph{disintegration} of $\mu$ with respect to  the Borel map $p : (Y, \mathcal{Y}) \to (Z, \mathcal{Z})$.

	Let $(\Omega, \mathcal{F}, \mathbb{P}, T, \varphi, X)$ be a discrete random dynamical system on the Polish space $(X, \mathcal{B})$, where we do not assume that $T$ is invertible. We let $\prob_\mathbb{P} (\Omega \times X )$ be the set of Borel probability measures $\mu$ on $\Omega \times X$ such that $(\pi_\Omega)_\ast\mu = \mathbb{P}$. Endow $\prob_\mathbb{P} (\Omega \times X )$ with the topology of weak convergence \cite[Definition~1.5.3]{arnold98}, with which $\prob_\mathbb{P} (\Omega \times X )$ is a Polish space since $(\Omega, \calF)$ is standard Borel \cite[Theorem 5.6]{crauel02}. 
	
	\begin{Def}
		Let $\mu \in \probP(\Omega \times X)$. Say that a function $\mu_\cdot (\cdot) : \Omega \times \mathcal{B} \to [0,1]$ is a \emph{disintegration} or \emph{factorization} of $\mu$ with respect to $\mathbb{P}$ if 
		\begin{enumerate}
			\item For all $Y \subseteq X$ Borel, $\omega \mapsto \mu_\omega (Y) $ is measurable. 
			\item For $\mathbb{P}$-almost every $\omega \in \Omega$, $\mu_\omega $ is a probability measure on $X$. 
			\item for all $A \in \mathcal{F} \otimes \mathcal{B}$, 
			$$ \mu(A) = \int_\Omega \int_X 1_A(\omega, x) d\mu_\omega(x) d\mathbb{P}(\omega).$$
			In other words, $d\mu(\omega, x) = d\mu_\omega (x) d\mathbb{P}(\omega)$
		\end{enumerate}
	\end{Def}

	We denote by $\Map (\Omega, \prob(X))$ the space of measurable maps $\mu_{\cdot} : \omega \in \Omega \mapsto \mu_\omega \in \prob(X)$, two such maps being identified if they agree almost everywhere. More generally, if $Y$ is a measured space and $Z$ is measurable, let $\Map(Y,Z)$ denote the set of equivalence classes of measurable maps from $X $ to $Y$, two such maps being identified if they agree almost everywhere. 
	If moreover $Y$ and $Z$ are $G$-spaces, we denote by $\Map_G(Y,Z) $ the set of $G$-equivariant measurable maps from $Y$ to $Z$, where two of them are identified if they agree almost everywhere. In case $X$ is compact, the Banach space $L^\infty (\Omega,M(X))$ with the variation norm on the space $M(X)$ of finite signed measures on $X$ is the dual of $L^1(\Omega, C(X))$ (with the supremum norm on the continuous functions $C(M)$). In this case, $\Map (\Omega, \prob(X))$ is a closed convex subset of $L^\infty (\Omega, M(X))$, hence compact.

	Note that if $X$ is a standard Borel space, any $\mu \in \probP(\Omega \times X)$ admits a disintegration  with respect to $\mathbb{P}$ \cite[Proposition 1.4.3]{arnold98}. 
	Conversely, any map  $\mu_\cdot \in \Map (\Omega, \prob(X))$ corresponds to a measure $\mu \in \probP(\Omega \times X)$ defined by 
	$$\mu(F \times B) = \int_F \mu_\omega (B) d\bbP(\omega), $$
	for any $F \times B \in \mathcal{F}\otimes \mathcal{B}$.

	If $\mathcal{E} \subseteq \mathcal{F} $ is a sub-$\sigma$-algebra, we can apply Rohlin's disintegration theorem for the factor $(\Omega \times X, \calF \otimes \calB) \to (\Omega \times X, \mathcal{E} \otimes \calB)$. If $\mu_\cdot \in L^\infty (\Omega, \calF; \prob(X))$ denote by $\overline{\mu}_\cdot \in L^\infty (\Omega, \mathcal{E}; \prob(X))$ the disintegration of the restriction of $\mu$ with respect to $\mathcal{E} \otimes \mathcal{B}$. For any $B \subseteq X$ measurable, we have that 
	$$ \overline{\mu}_\omega (B) = \mathbb{E}(\mu_\cdot (B) | \mathcal{E})(\omega), $$
	where we denoted by $\mathbb{E}(\mu_\cdot (B) | \mathcal{E})$ the conditional expectation of the random variable $\omega \mapsto \mu_\omega(B)$. 
	It is the \emph{conditional expectation of $\mu$} with respect to $\mathcal{E} \otimes \mathcal{B}$ and is denoted by $\omega \mapsto \mathbb{E}(\mu_\cdot | \mathcal{E})_\omega$. We have that
	\begin{align}
		\mathbb{E}(\mu_\cdot (B) | \mathcal{E})(\omega) = \mathbb{E}(\mu_\cdot  | \mathcal{E})_\omega (B). \nonumber
	\end{align}

\subsection{Invariant measures}
	Let $(\Omega, \mathcal{F}, \mathbb{P}, T, \varphi, X)$ be a RDS on the Polish space $(X, \mathcal{B})$. Consider the skew product 
	\begin{align}
		T^X \colon & \Omega \times X \to \Omega \times X \nonumber \\
		& (\omega, x) \mapsto (T\omega, \varphi(1, \omega) x) \nonumber.
	\end{align}
	
	\begin{Def}
		We say that a measure $\nu \in \prob_\mathbb{P} (\Omega \times X )$ is \emph{$T^X$-invariant} if $(T^X)_\ast \nu = \nu$. The set of $T^X$-invariant probability measures is denoted by $\mathcal{I}_\mathbb{P}(\varphi)$. 
	\end{Def}
	Understanding invariant measures for the RDS is a crucial step in studying its asymptotic behavior. Fortunately, when the Polish space $(X, \calB)$ is compact, there always exist invariant measures. 
	\begin{thm}[{Krylov-Bogolyubov Procedure, \cite[Theorem 1.5.8]{arnold98}}]\label{thm krylov}
		Let $\varphi$ be a continuous RDS over $(\Omega, \mathcal{F}, \mathbb{P}, T,X)$, i.e., a RDS such that for all $n, \omega$, $\varphi(n, \omega ) $ is a continuous mapping. Fix a probability measure $m \in \probP(\Omega \times X)$. Let
		\begin{align}
			m_N := \frac{1}{N} \sum_{n=0}^{N-1} ((T^X)^n)_\ast m. \nonumber
		\end{align}
		Then any limit point of $m_N$ as $N \rightarrow \infty$ in the topology of weak-$\ast$ convergence is a $T^X$-invariant probability measure in $\probP(\Omega \times X)$. 
	\end{thm}
	
	\begin{cor}
		If $X$ is moreover compact, there exists an invariant probability measure $\nu \in \invmeas$. 
	\end{cor}
	
	\begin{proof}
		Let $\eta$ be any Borel probability measure on $X$. Apply the Krylov-Bogolyubov Procedure to the probability measure $m = \mathbb{P}\otimes \eta$. As $X$ is compact, $\Map (\Omega,\prob(X))$ is convex compact and there exists a subsequence of $m_N$ which is convergent for the weak-$\ast$ topology. 
	\end{proof}
	
	The following lemma, whose proof is straightforward, relates invariance and factorization \cite[Theorem 1.4.5]{arnold98}. 
	
	\begin{lem}\label{lem inv fact}
		A measure $\mu \in \probP(\Omega \times X)$ is $T^X$-invariant if and only if for all $n \in \mathbb{N}$,
		$$\mathbb{E}(\varphi(n, \cdot)\mu_\cdot | (T^{n})^{-1}\mathcal{F}) (\omega) = \mu_{T^n\omega},$$ 
		where $\mu_\cdot \in L^\infty (\Omega, \calF; \prob(X))$ is the disintegration associated to $\mu$ with respect to $\bbP$. 
		
		In particular, if the RDS $\varphi$ is invertible, $T^{-1} \calF = \calF$ and $\mu$ is $T^X$-invariant if and only if for all $n \in \mathbb{Z}$
		$$\varphi(n, \omega)_\ast\mu_\omega = \mu_{T^n\omega}.$$
	\end{lem}
	
	\begin{rem}
		If the flow $T$ is discrete, the conditions in Lemma \ref{lem inv fact} are equivalent to the simpler 
		$$ \mathbb{E}(\varphi(1, \cdot)\mu_\cdot | T^{-1}\mathcal{F}) (\omega) = \mu_{T\omega} \ \text{ $\bbP$-almost surely}, $$
		and in case $T$ is invertible, to 
		$$\varphi(1, \omega)_\ast \mu_\omega = \mu_{T\omega}  \ \text{ $\bbP$-almost surely}.$$
	\end{rem}

\subsection{The natural extension}
	 Consider the standard probability space $(\Omega, \mathcal{F}, \mathbb{P})$, and let $T: \Omega \to \Omega$ be a measure-preserving transformation, not necessarily invertible. 
	\begin{Def}
		The \emph{natural extension} of $(\Omega, \mathcal{F}, \mathbb{P}, T)$ is the system $(\tilde{\Omega}, \tilde{\calF}, \tilde{\bbP}, \tilde{T})$ where 
		\begin{enumerate}
			\item $\tilde{\Omega}:= \{\bar{\omega} = (\omega_i)_{i \in \bbZ} \mid  \omega_i \in \Omega , T\omega_i = \omega_{i+1} \ \forall i\}$;
			\item $\tilde{\calF}$ is the smallest $\sigma$-algebra which contains all sets of the form $\{\bar{\omega} \in \tilde{\Omega} \mid \omega_i \in E\}$, with $i \leq 0 $ and $E \in \calB$;
			\item $\tilde{\bbP}$ is the unique probability measure on $\tilde{\calF}$ such that $\tilde{\bbP}\{\bar{x} \in \tilde{X} \mid x_i \in E_i\} = \mu(E_i)$ for all $i \leq 0 $ and $E_i \in T^{-i}\calB$; 
			\item $\tilde{T}$ is the left Bernoulli shift: for all $\bar{\omega} = (\omega_i)\in \tilde{X}$, $T\bar{\omega} = (\omega_{i+1})_{i \in \bbZ}$. 
		\end{enumerate}
	\end{Def}
	The following result is standard in ergodic theory, see for instance \cite[\S1.6.3]{sarig23}. 
	\begin{thm}
		Consider the same setting and assumptions. Then the probability measure $\tilde{\bbP}$ exits and is unique. The natural extension $(\tilde{\Omega}, \tilde{\calF}, \tilde{\bbP}, \tilde{T})$ of $(X, \calB, \mu, T)$ is an invertible extension of $(\Omega, \mathcal{F}, \mathbb{P}, T)$, and is the factor of any other invertible extension of $(\Omega, \mathcal{F}, \mathbb{P}, T)$. The probability measure preserving transformation $\tilde{T}$ is ergodic if and only if $T$ is ergodic, and mixing if and only if $T$ is mixing. 
	\end{thm}
	
	In particular, starting from a probability-preserving transformation on a standard Borel space $(\Omega, \mathcal{F}, \mathbb{P}, T)$, we can always consider the natural invertible extension. 
	
	\subsection{Past and future factors}\label{section past future sigma algebras}
	
	In this section we do the converse of the natural extension in that we reduce bilateral processes to unilateral ones: we describe how to restrict an invertible RDS to positive or negative times. Consider a discrete RDS $(\Omega, \mathcal{F}, \mathbb{P}, T, \varphi, X)$, where $T$ is invertible (in our case, it means that $\mathbb{T}= \bbZ$), and $(X, \calB)$ is a Polish space with its $\sigma$-algebra of Borel sets and we assume that the RDS $\varphi$ has values in the group of homeomorphisms of $X$. 
	
	This datum defines the following two-parameter filtration on $(\Omega, \calF, \bbP)$: for all $-\infty\leq s\leq t\leq +\infty  $
	$$ \calF_s^t = \sigma\{{\omega\mapsto \varphi(\tau, T^{s }\omega)x \mid   0\leq \tau\leq t-s\, ; \ x \in X}\}.$$
	
	In particular, the $\sigma$-algebra $\mathcal{F}_0:= \calF_0^1$ is the smallest $\mathbb{P}$-complete (i.e. all null sets are measurable) $\sigma$-algebra on which $f= \varphi(1, \cdot)$ is measurable. Recall that $T$ acts on the set of $\mathbb{P}$-complete $\sigma$-subalgebras $\mathcal{F}'$ of $\mathcal{F}$ by $$T^{-1}\mathcal{F}' := \{T^{-1}Y\mid Y \in \mathcal{F}'\}.$$
	This definition also holds whenever $T$ is just assumed measure-preserving, but not necessarily invertible. Likewise, $\mathcal{F}_n^m$ is the smallest $\mathbb{P}$-complete $\sigma$-algebra on which $f \circ T^{n}, \dots , f \circ T^{m}$ is measurable. In particular, for $n \geq 1$, $\varphi(n, .)$ is $\mathcal{F}_0^{n}$-measurable. 
	We can assume without loss of generality that $\mathcal{F}_{-\infty}^{+\infty} = \mathcal{F}$. Consider the sub-$\sigma$-algebras $\calF^+:= \mathcal{F}_{0}^{+\infty}$ and $\calF^-:= \mathcal{F}_{-\infty}^{-1}$. The $\sigma$-algebra $\calF^+$ is $\bbP$-complete, and is such that for all $n \geq 1$,  
	\begin{enumerate}
		\item $T^{-n} \calF^+ \subseteq \calF^+$;
		\item $\omega \mapsto \varphi(n, \omega)$ is $\calF^+$-measurable. 
	\end{enumerate}
	Hence we can consider the measurable \emph{future factor} associated 
	\begin{align}
		p_+ : (\Omega, \mathcal{F})\to (\Omega, \mathcal{F}^+),  \nonumber
	\end{align}
	and restrict the two-sided RDS $(\Omega, \calF, T, \varphi)$ to the one-sided RDS $(\Omega, \bbP_+, \calF^+, T_+, \varphi)$. 
	
	The analogue is true for $\calF^-$, with the factor 
	$$p_- : (\Omega, \mathcal{F})\to (\Omega, \mathcal{F}_-).$$
	In order to differentiate with $(\Omega, \mathcal{F})$, we denote by $\Omega_+$ and $\Omega_-$ the standard Borel spaces $(\Omega, \calF^+, \mathbb{P}^+ := (p_{+})_\ast\mathbb{P})$ and $(\Omega, \calF^-, \mathbb{P}^- := (p_{-})_\ast\mathbb{P})$ with these measures respectively. 
	
	By construction, the transformation $T_+ $ satisfies 
	$$T_+(p_+(\omega)) = p_+(T\omega)$$
	almost surely. We call $(\Omega, \mathbb{P}^+, \calF^+)$ the \emph{future factor} of the greg $(\Omega, \mathcal{F}, \mathbb{P}, T, \varphi)$, and the one-sided RDS associated is the \emph{forward RDS}. 
	Similarly, the probability space $(\Omega_-, \mathbb{P}^-, \calF^-)$ is called the \emph{past factor} associated to $\mathcal{F}_{-\infty}^{0}$ and satisfies 
	$$T_-(p_-(\omega)) = p_-(T^{-1}\omega)$$
	almost surely. The one-sided RDS associated to $(\Omega_-, \mathbb{P}^-, \calF^-, T_-, f)$ is called the \emph{backward RDS}.

	\subsection{Forward invariant measures and martingales}
	
	Let $(\Omega, \mathbb{P}, T, G, f)$ be an invertible discrete RDS on $G$ and let $G$ act on the Polish space $(X, \calB)$. 
	We define
	\begin{align}
		\Map_f(\Omega, \prob{X}) = \{ \mu_\cdot \in \Map(\Omega, \prob{X}) \mid  \mu_{T \omega}= f(\omega)\mu_{\omega} \text{ almost everywhere}\} \nonumber.
	\end{align}
	Lemma \ref{lem inv fact} states that $\Map_f(\Omega, X) $ is exactly the set of measures $\mu \in \probP(\Omega \times X)$ that are  $T^X$-invariant.

	As explained in the previous section, we can restrict the RDS to non-negative times by considering the sub-$\sigma$-algebra $\calF^+$. 
	Let us denote by $\Map (\Omega, \calF^+; \prob(X))$ the space 
	$$ \{\mu^+_\cdot \in \Map (\Omega,\prob(X)) \mid \mu^+_\cdot \text{ is measurable w.r.t. } \calF^+\}.$$
	
	A measure $\mu^+_\cdot \in \Map(\Omega, \calF^+; \prob(X))$ is \emph{forward invariant} if it is invariant for the non-invertible RDS $(\Omega, \calF^+, \mathbb{P}^+, T_+, \varphi)$.
	
	By Lemma \ref{lem inv fact}, a forward invariant measure $\mu^+_\cdot$ satisfies: for $\mathbb{P}^+$-almost every $\omega \in \Omega_+$, 
	$$\mathbb{E}(f(\cdot)_\ast\mu^+_\cdot| T^{-1}\mathcal{F}^+)(\omega) = \mu^+_{T\omega}.$$

	\begin{lem}\label{lem eq stat forward invariant invertible}
		Suppose $T$ is invertible. A map $\mu^+_\cdot \in \Map (\Omega, \calF^+; \prob(X))$ is forward-invariant if and only if for all $n \geq 1$,
		\begin{align}\label{eq stat forward invariant invertible}
			\mu^+_\omega = \mathbb{E}(\varphi(-n, \cdot)^{-1}\mu^+_{T^{-n}\cdot} | \mathcal{F}^+)(\omega). 
		\end{align}
	\end{lem}
	\begin{proof}
		As $\mu^+_\cdot \in \Map (\Omega; \prob(X))$, we can use that $T$ is invertible and we obtain the following equivalent characterization: 
		Indeed, by disintegration, $\mu^+_\cdot$ can be seen as a $\calF^+$-measurable map in $\Map (\Omega, \prob(X))$, and it is forward-invariant if and only almost surely, for all $n \in \bbN$, 
		$$\mathbb{E}(\varphi(n, \cdot)\mu^+_\cdot | (T^{n})^{-1}\mathcal{F}^+) (\omega) = \mu^+_{T^n\omega},$$
		hence \begin{align}
			\mu^+_{\omega} &= (\mu^+_{T^n\cdot} \circ T^{-n}) (\omega) = \mathbb{E}(\varphi(n, \cdot)\mu^+_\cdot | (T^{n})^{-1}\mathcal{F}^+) \circ T^{-n}(\omega) \nonumber \\
			& = \mathbb{E}((\varphi(n, \cdot)\mu^+_\cdot)\circ T^{-n} | \mathcal{F}^+)(\omega)\nonumber \\
			& =  \mathbb{E}(\varphi(n, T^{-n}\cdot)\mu^+_{T^{-n}\cdot} | \mathcal{F}^+) (\omega)\nonumber,
		\end{align}
		and $\varphi(n, T^{-n}\cdot) = \varphi(-n, \cdot)^{-1}$ almost surely by the cocycle property. 
	\end{proof}

	Let us rewrite the preceding equation. Denote by $\omega \mapsto \mathbb{P}^+_\omega $ the disintegration of $\mathbb{P}$ associated to the Borel factor $p_+$. A $\mathcal{F}^+$-measurable map $\mu^+_\cdot : \Omega \to \prob(X)$ is stationary if for any non-negative measurable function $h: X \to \bbR_+$, for almost every $\omega \in \Omega$, 
	$$\int_X h (x) d\mu^+_\omega(x) = \int_\Omega \int_X h(f(T^{-1}\omega') x ) d\mu^+_{T^{-1}\omega}(x) d\mathbb{P}^+_\omega(\omega'). $$
	
	In particular, if $\mu^+_\cdot$ is forward-invariant then for all $U\in \calB$, for all $n \geq 1$
	$$\mu^+_\omega (U) = \int_\Omega \mu^+_{T^{-n}\omega'} (\varphi(n, T^{-n}\omega')^{-1}U) d\bbP^+_\omega (\omega')$$ 
	for $\bbP$-almost every $\omega$. 
	We call these maps \emph{forward stationary maps}, and we denote the space of forward stationary maps by $\Mapstat(\Omega_+, \prob(X))$.

	The following theorem is the RDS analogue of a Theorem due to Furstenberg about limit measures for random walks, and will be crucial in the sequel. In the context of dynamical system, it appears in \cite[Lemma 1]{ledrappier86} for products of independent and identically distributed diffeomorphisms of a projective space. Its proof also relies on Doob's Martingale Convergence Theorem.  
	\begin{thm}\label{thm lim measures rds}
		There is a bijection between invariant measures $\mu \in \invmeas$ and forward invariant measures $\mathcal{I}_\mathbb{P}(\varphi^+)$. More precisely, for $\mu \in \invmeas$, the measure $\mu^+_\cdot \in \probP(\Omega_+ \times X)$ defined by
		$$ \mu^+_\omega = \mathbb{E}(\mu_\cdot | \calF^+)(\omega)$$
		is forward invariant. Conversely, for $\mu^+_\cdot \in L^\infty (\Omega, \calF^+; \prob(X))$ forward invariant, we have the a.e. weak-$\ast$ convergence
		\begin{align}
			\mu_\omega = \lim_n \varphi(n, T^{-n}\omega)\mu^+_{p^+(T^{-n}\omega)}\nonumber,
		\end{align}	
		and the measure $\mu_\cdot \in L^\infty (\Omega; \prob(X))$ is invariant for the two-sided RDS. 
	\end{thm}
	
	\begin{proof}
		Let $\mu$ be an invariant measure for the RDS $(\Omega, \mathcal{F}, \mathbb{P}, T,\varphi, X)$. Consider the measure $\mu^+ = \mu\mid_{\calB \otimes \calF^+}$. If we denote by $\omega \mapsto \mu_\omega^+$ the disintegration with respect to the restriction of $\calF$ to the $\sigma$-algebra $\calF^+ $, we have as before 
		$$ \mu^+_\omega = \mathbb{E}(\mu_\cdot | \calF^+)(\omega).$$
		We have that for almost every $\omega \in \Omega$ and for all $n \in \bbN$, 
		\begin{align}
			\mathbb{E}(\varphi(-n, \cdot)^{-1}\mu^+_{T^{-n}\cdot} | \mathcal{F}^+)_\omega & = \mathbb{E}(\varphi(-n, \cdot)^{-1}\mathbb{E}(\mu_{T^{-n}\cdot}|\calF^+) | \mathcal{F}^+)_\omega\nonumber \\
			& =  \mathbb{E}(\mathbb{E}(\varphi(-n, \cdot)^{-1}\mu_{T^{-n}\cdot}|\calF^+) | \mathcal{F}^+)_\omega \nonumber \\
			& =  \mathbb{E}(\mu_{\cdot}| \mathcal{F}^+)_\omega \text{ by Lemma \ref{lem inv fact}}\nonumber \\
			& =    \mu^+_\omega \nonumber. 
		\end{align}
		By Lemma \ref{lem eq stat forward invariant invertible}, $\mu^+$ is forward-invariant. 
		
		Conversely, let $\mu^+ \in L^\infty (\Omega, \calF^+; \prob(X))$ be a forward-invariant measure, which we see as maps in $L^\infty (\Omega; \prob(X))$ that is $\calF^+$-invariant (and we write $\mu^+_{p_+(\omega)} = \mu^+_\omega$ to simplify the notations). Recall that by construction, $(T^{n})^{-1}\calF^+ = \calF_n^\infty$. We are going to show that seen as maps in $L^\infty (\Omega; \prob(X))$, the sequence $(\chi(n, \cdot)\mu^+_{T^{-n}\cdot})$ forms a martingale with respect to the filtration $\calF^\infty_{-n}$. Indeed, 
		\begin{align}
			& \mathbb{E}(\chi(n+1, \cdot)\mu^+_{T^{-(n+1)}\cdot}| \calF^\infty_{-n})_\omega = \mathbb{E}(\chi(n+1, \cdot)\mu^+_{T^{-(n+1)}\cdot}| \calF^\infty_{-n})_\omega \nonumber \\
			& =   \mathbb{E}(\chi(n, \cdot)\chi(1, T^{-n}\cdot)\mu^+_{T^{-(n+1)}\cdot}| \calF^\infty_{-n})_\omega \nonumber \\
			& =  \chi(n, \omega) \mathbb{E}(f(T^{-(n+1)}\cdot)\mu^+_{T^{-(n+1)}\cdot}| \calF^\infty_{-n})_\omega \text{ because $\chi(n, \cdot)$ is $\calF^\infty_{-n}$-measurable}, \nonumber \\
			& = \chi(n, \omega) \mathbb{E}((f(\cdot)\mu^+_{\cdot})\circ T^{-(n+1)}| \calF^\infty_{-n})_\omega \nonumber \\
			& =  \chi(n, \omega) \mathbb{E}(f(\cdot)\mu^+_{\cdot}| \calF^\infty_{1}) \circ T^{-(n+1)}(\omega )\nonumber \\
			& =  \chi(n, \omega) (\mu^+_{T\cdot}\circ T^{-(n+1)})(\omega )\text{ by Lemma \ref{lem inv fact}} \nonumber \\
			& =  \chi(n,\omega) \mu^+_{T^{-n}\omega} \nonumber.
		\end{align}
		By Doobs martingale convergence theorem, the limit 
		\begin{align}
			\mu_\omega := \lim_{n \to \infty}\chi(n,\omega) \mu^+_{T^{-n}\omega} \label{eq limit forward}
		\end{align} 
		exists almost surely. This martingale is uniformly integrable, so the convergence takes place in $L^1(\bbP)$. This yields a measurable map $\mu_\cdot \in L^\infty (\Omega, \calF; \prob(X))$, where $\calF = \calF_{-\infty}^{+\infty}$. Moreover, for all $n$, 
		\begin{align}
			\mathbb{E}(\mu_\cdot \mid \calF_{-n}^\infty) = \chi(n, \cdot) \mu^+_{T^{-n}\cdot}. \label{eq bijection conditional}
		\end{align}
		
		The measure $\mu$ is invariant by $\varphi$-invariance of the limit. Indeed, as for all $m$, the map 
		$$ \Theta(m) : \nu_\cdot \in \prob_\bbP(\Omega \times X) \mapsto \chi(m, \cdot) \nu_{T^{-m}\cdot} \in \prob_\bbP(\Omega \times X)$$ 
		is continuous by continuity of the maps $\chi(m, \cdot)$. Hence for all $m \in \bbZ$, the limit  \eqref{eq limit forward} gives 
		$$\mu_\cdot = \lim \Theta(n) \mu^+_\cdot = \lim_{n\to \infty} \Theta(m) \Theta(n) \mu^+_\cdot = \Theta(m)  \lim_{n\to \infty} \Theta(n) \mu^+_\cdot = \Theta(m) \mu_\cdot \, $$ 
		showing that $\mu_\cdot$ is indeed invariant. 
		
		As $\calF = \calF_{-\infty}^{+\infty}$, Equation \ref{eq bijection conditional} shows that if $\mu \in \prob_\bbP(\Omega \times X)$ is invariant, the procedure $\mu_\cdot \in  \mapsto \mathbb{E}(\mu_\cdot\mid \calF) \mapsto \lim \chi(n, \cdot) \mu^+_{T^{-n}\cdot}$ is the identity, showing the bijection. 
	\end{proof}

	\section{Boundary systems and Mackey range of a cocycle}\label{section bd system}

	Let $(\Omega, \mathcal{F}, \mathbb{P}, T, f, G)$ be a greg on a locally compact second countable group $G$. In this section, we construct the space of ideal pasts $B_-$ and the space of ideal futures $B_+$ associated to the RDS $\varphi$: they are the spaces of ergodic components for the natural $T_-$-action on $\Omega_-$ and for the natural $T_+$-action on $\Omega_+$ on the past and future factors respectively. These spaces are the cocycle analogues of the Poisson boundary $B(G, \mu)$ associated to the random walk generated by an admissible measure $\mu$ on the group $G$. 
	
	\subsection{Spaces of trajectories}\label{section space trajectories}
	We assume the group $G$ to be locally compact and second countable. As a consequence, we can (and shall) take a probability measure $\mu_G$ on $G $ in the measure class of the right Haar measure. Let us consider the map 
	\begin{align}
		\alpha \colon & \Omega \times G \to G^\mathbb{Z}\nonumber \\
		& (\omega,g) \mapsto \alpha_g(\omega)=(\varphi(n, \omega)g^{-1})_{n \in \mathbb{Z} }. \nonumber
	\end{align}

	Endow the target space with the Lebesgue measure $\tilde{\mathbb{P}} = \alpha_\ast ( \mathbb{P}\otimes \mu_G)$, and write $\tilde{\Omega} =  G^\mathbb{Z} $ as the measure space with the measure $\tilde{\mathbb{P}}$. The Lebesgue space $\tilde{\Omega}$ is the \emph{space of trajectories} of the ergodic cocycle $(\varphi(n, \omega)g^{-1})_{n\in \mathbb{Z}}$, with initial value $g^{-1} \in G$ distributed according to $\mu_G$. Similarly, denote by $\tilde{\mathbb{P}}_g = (\alpha_{g})_\ast \mathbb{P}$, the image of $\mathbb{P}$ under the measurable map 
	\begin{align}
		\alpha_g \colon & \Omega \to \tilde{\Omega}\nonumber\\
		& \omega \mapsto (\varphi(n, \omega)g^{-1})_{n \in \mathbb{Z} }. \nonumber
	\end{align} 
	Note that $\tilde{\mathbb{P}}_g$ is a probability measure on the space of trajectories $\tilde{\Omega}$, and that 
	$$\tilde{\mathbb{P}} = \int_G \tilde{\mathbb{P}}_g d\mu_G(g).$$
	Since we assumed that $\mathcal{F}_{-\infty}^\infty = \mathcal{F}$, the measurable map 
	\begin{align}
		\alpha_e \colon& \Omega \to \tilde{\Omega}\nonumber\\
		& \omega \mapsto (\varphi(n, \omega))_{n \in \mathbb{Z} } \nonumber
	\end{align}
	is almost surely injective. In particular, $(\Omega, \mathbb{P})$ and $(\tilde{\Omega},\tilde{\mathbb{P}}_e)$ are isomorphic as standard Borel spaces. The same is true between $(\Omega, \mathbb{P})$ and $(\tilde{\Omega},\tilde{\mathbb{P}}_g)$ for any $g$.
	
	The natural $G$-action on $\tilde{\Omega} $ is given by 
	\begin{align}
		 & G \times \tilde{\Omega} \rightarrow \tilde{\Omega} \nonumber\\
		 &(g, (g_n)_{n\in\mathbb{Z}}) \mapsto  (g_ng^{-1})_{n\in\mathbb{Z}} \nonumber, 
	\end{align}
	so that with the natural $G$-action on $\Omega \times G$ given by $g \cdot (\omega, h) = (\omega, gh)$, the map $\alpha$ is $G$-equivariant. The $\mathbb{Z}$-action on $\tilde{\Omega}$ is given by the shift-map 
	\begin{align}
		S \colon&\tilde{\Omega} \rightarrow \tilde{\Omega} \nonumber\\
		& (g_n)_{n\in\mathbb{Z}} \mapsto (g_{n+1})_{n\in\mathbb{Z}}. \nonumber
	\end{align}
	
	These actions commute and, because $\mu_G$ is equivalent to the Haar measure on $G$, they are non-singular, i.e measure-class preserving. Moreover, the flow $T$ and the shift map $S$ are intertwined in the following way: for $\mu_G$-almost every $g$ and $\mathbb{P}$-almost every $\omega \in\Omega$, 
	\begin{align}
		\alpha(g, T\omega) &= (\varphi(n, T\omega)g^{-1})_{n \in \mathbb{Z} } \nonumber \\
		& =  (\varphi(n+1, \omega)\varphi(1, \omega)^{-1}g^{-1})_{n \in \mathbb{Z} } \nonumber \\
		& =  (\varphi(n+1, \omega)(gf(\omega))^{-1})_{n \in \mathbb{Z} } \nonumber \\
		& =  S\alpha(gf(\omega), \omega). \nonumber
	\end{align}
	As a consequence, for $\mu_G$-almost every $g$, $\mathbb{P}$-almost every $\omega \in\Omega$ and all $m \in \mathbb{Z}$, 
	\begin{align}
		\alpha_g (T^{m}\omega) = S^m\alpha_{g\varphi(m,\omega)}( \omega), \nonumber
	\end{align}
	or in the other direction 
	\begin{align}
		S^m \alpha_g (\omega) = \alpha_{g\varphi(m,\omega)^{-1}}( T^m \omega).\nonumber
	\end{align}
	We can repeat the same constructions with the future and past factors. Indeed, consider the map 
	\begin{align}
		\alpha^+ \colon& \Omega_+ \times G \to G^\mathbb{N}\nonumber \\
		& (\omega,g) \mapsto \alpha_g(\omega)=(\varphi(n, \omega)g^{-1})_{n \in \mathbb{N} }. \nonumber
	\end{align}
	and endow the target space $\tilde{\Omega}_+ := G^\mathbb{N}$ with the measure $\tilde{\mathbb{P}}_+ = (\alpha^+)_\ast ( \mathbb{P}^+\otimes \mu_G)$. The right $G$-action on $\tilde{\Omega}_+$ remains the same and when we consider the restriction $S_+: \tilde{\Omega}_+ \to \tilde{\Omega}_+$ of the shift map to positive times, so that for $m \geq 0$,
	\begin{align}
		\alpha^+_g ((T_+)^{m}\omega) = (S_+)^m\alpha^+_{g\varphi(m,\omega)}( \omega).\nonumber
	\end{align}
	Similarly, $\tilde{\Omega}_- := G^\mathbb{N}$ is the Lebesgue space with the measure $\tilde{\mathbb{P}}_- = (\alpha^-)_\ast ( \mathbb{P}^-\otimes \mu_G)$, where
	\begin{align}
		\alpha^- \colon & \Omega_- \times G \to \tilde{\Omega}_- \nonumber \\
		& (\omega,g) \mapsto \alpha^-_g(\omega)=(\varphi(n, \omega)g^{-1})_{n \leq -1 }. \nonumber
	\end{align}
	Define the inverse shift map $S_-: \tilde{\Omega}_- \to \tilde{\Omega}_-$ as $(g_i)_{i\in\mathbb{N}} \mapsto (g_{i-1})_{i\in\mathbb{N}} $. 
	
	We have the intertwining: for $m \geq 0 $, 
	\begin{align}
		\alpha^-_g ((T_-)^{m}\omega) = (S_-)^m\alpha^-_{g\varphi(-m,\omega)}( \omega). \nonumber
	\end{align}
	or in the other direction 
	\begin{align}
		(S_-)^m \alpha^-_g (\omega) = \alpha^-_{g\varphi(-m,\omega)^{-1}}( (T_-)^{m} \omega).\nonumber
	\end{align}
	
	Notice that the projections $(\tilde{\Omega}, \tilde{P}) \to (\tilde{\Omega}_+, \tilde{P}_+)$ and $(\tilde{\Omega}, \tilde{P}) \to (\tilde{\Omega}_-, \tilde{P}_-)$ are $G$-equivariant, measure-preserving and $S_+$ or $S_-$-equivariant respectively. 
	\newline
	
	For the sake of clarity in Section~\ref{section erg prop mackey range}, we shall also consider the space of trajectories associated with the backward cocycle $\chi $. Namely, define 
	\begin{align}
		\check{\alpha} \colon & \Omega \times G \to G^\mathbb{Z}\nonumber \\
		& (\omega,g) \mapsto \check{\alpha}_g(\omega)=(g\chi(n, \omega))_{n \in \mathbb{Z} }. \nonumber
	\end{align}
	Define also $(\check{\Omega}, \check{\mathbb{P}} ):= (G^\mathbb{Z},\check{\alpha}_\ast( \mathbb{P}\otimes \mu_G))$.  The natural $G$-action on $\check{\Omega} $ is given by 
	\begin{align}
		& G \times \check{\Omega} \rightarrow \check{\Omega} \nonumber\\
		&(g, (g_n)_{n\in\mathbb{Z}}) \mapsto  (g g_n)_{n\in\mathbb{Z}} \nonumber, 
	\end{align}
	while the $\mathbb{Z}$-action on $\check{\Omega}$ is given by the shift-map 
	\begin{align}
		\check{S} \colon&\check{\Omega} \rightarrow \check{\Omega} \nonumber\\
		& (g_n)_{n\in\mathbb{Z}} \mapsto (g_{n+1})_{n\in\mathbb{Z}}. \nonumber
	\end{align}
	These actions commute and are non-singular. Define $\check{T}:= T^{-1}$ so that $\chi $ is a backward $\check{T}$-cocycle: for $n, m \in \bbZ$ and $\omega \in \Omega$,
	\[\chi(n+m, \omega)= \chi(m, \omega) \chi (n, \Ti^m \omega ).\] We then have the intertwining
	\begin{align}
		\Si^m\alphai_g(\omega) = \alphai_{g\chi(m, \omega)} (\Ti^m \omega) \label{eq intertwining chi}.
	\end{align}
	Of course, we define $(\Omegai_+, \bbPi_+)$ and $(\Omegai_-,\bbPi_-)$ in a similar fashion as $(\Omega_+,\tilde{\mathbb{P}}_+) $ and $(\Omega_-, \tilde{\mathbb{P}}_-)$, with the maps
	\begin{align}
		\alphai^+ \colon& \Omega_- \times G \to \Omegai_+ \nonumber \\
		& (\omega,g) \mapsto \alphai^+_g(\omega)=(g \chi(n, \omega))_{n \in \mathbb{N} } \nonumber
	\end{align}
	and 
	\begin{align}
		\alphai^- \colon& \Omega_+ \times G \to \Omegai_-\nonumber \\
		& (\omega,g) \mapsto \alphai^-_g(\omega)=(g \chi(n, \omega))_{n\leq -1}. \nonumber
	\end{align}
	In particular, the spaces $(\Omega_-, \bbP_-)$ and $ (\Omegai_+, \bbPi_e)$ are as standard Borel spaces. 
	We finish by the following basic observation. 
	\begin{lem}\label{lem isom traj chi -}
		The transformation 
		\begin{align}
			\iota\colon& \Omegai_+ \rightarrow \tilde{\Omega}_-  \nonumber\\
			& (g_n)_{n\in\mathbb{N}} \mapsto (g_n^{-1})_{n\in\mathbb{N}}. \nonumber
		\end{align}
		is an isomorphism between the standard measured spaces $(\Omega_-, \tilde{\mathbb{P}}_-)$ and $(\Omegai_+, \bbPi_+)$, and satisfies $\iota(\Si \omegai) = S_-(\iota (\omegai))$ for $\omegai \in \Omegai$. 
		An analogous statement holds for the spaces $\Omegai_- $ and $\tilde{\Omega}_+$. 
	\end{lem}
	\begin{proof}
		For $n \geq 0 $ and $\omega \in \Omega$, we have by definition
		\[ \chi (n, \omega) = \varphi(-n, \omega)^{-1}, \] which proves the claim. 
	\end{proof}
	\subsection{Ideal futures and ideal pasts}
	We are now going to build the equivalent of the Poisson-Furstenberg boundary in the context of ergodic cocycles. The following construction relies on ideas of Mackey \cite{mackey62} and Zimmer \cite[\S5]{zimmer78}.
	
	Consider the set $$L^\infty(\tilde{\Omega}_+, \tilde{\bbP}_{+})^{S_{+}} := \{F \in L^\infty(\tilde{\Omega}_{+}, \tilde{\bbP})\mid F\circ S_{+} = F\}.$$ 
	Since the $G$-actions and $S_+$-actions on $\tilde{\Omega}_{+}$ commute, $L^\infty(\tilde{\Omega}_{+}, \tilde{\bbP}_{+})^{S_{+}}$ is a $G$-invariant von Neumann subalgebra. The following proposition is Mackey's point realization Theorem \cite[Section 3]{mackey62}, see also \cite[Appendix B]{zimmer84} or \cite[\S1]{takesaki03c}. 
	
	\begin{prop}[{Mackey}]
		Let $G \curvearrowright (X, \nu)$ be a non-singular action on a standard probability space $(X, \nu)$. Let $A \subseteq L^\infty (X, \nu)$ be a $G$-invariant von Neumann subalgebra. 
		
		Then there exist a compact metrizable space $Z$, a continuous action $G \curvearrowright Z$, a measure $\eta \in \prob(Z)$ and a $G$-equivariant measurable factor map $\pi : (X, \nu) \to (Z, \eta)$ so that the $G$-equivariant weak-$\ast$-continuous unital $\ast$-homomorphism $\pi^\ast : L^\infty (Z, \eta) \to L^\infty (X, \nu) $ given by $\pi^\ast(f) = f \circ \pi $ satisfies $\nu \circ \pi^\ast = \eta $ and
		$$\pi^\ast(L^\infty (Z) ) = A.$$
	\end{prop}
	
	As a consequence, there exists a standard Borel space $(B_+, \nu_+) $ with a distinguished measure class and endowed with a non-singular $G$-action such that $B_+$ realizes the quotient space of $\tilde{\Omega}_+$ by the partition given by $S_+$-stationary measurable subsets. This space is essentially unique \cite[Section 5]{mackey62}, and there exists a canonical $G$-equivariant measurable map 
	\begin{align}
		\bnd_+ : \tilde{\Omega}_+ \rightarrow B_+ .\nonumber
	\end{align}
	
	The space $(B_+, \nu_+) $ describes the (future) asymptotic behaviour of the ergodic cocycle $(\varphi(n, \omega))_n$. It is the space of ergodic components of $\tilde{\Omega}_+$ under the action of the semigroup $\mathbb{N}$ with $S_+$, and is sometimes written 
	$$ (B_+, \nu_+) = (\tilde{\Omega}_+, \tilde{\mathbb{P}}_+)\sslash S_+.$$
	We can build $(B_-, \nu_-)$ in a similar fashion by considering $L^\infty(\tilde{\Omega}_-, \tilde{\bbP}_{-})^{S_{-}}$. we can then define:
	$$ (B_-, \nu_-) = (\tilde{\Omega}_-, \tilde{\mathbb{P}}_-)\sslash S_-.$$
	Finally, there exists a Lebesgue probability space $(B, \nu) $ endowed with a non-singular $G$-action such that $B$ realizes the ergodic components of $\tilde{\Omega}\sslash S$: 
	$$ (B, \nu) = (\tilde{\Omega}, \tilde{\mathbb{P}})\sslash S.$$
	
	As the construction is functorial and there exist $S$-equivariant maps $\tilde{\Omega} \to \tilde{\Omega}_\pm$, We have natural $G$-equivariant measurable mappings
	\begin{align}
		\pi_+ : B \rightarrow B_+ \emph{ and } \pi_- : B \rightarrow B_-. \nonumber
	\end{align}

	We define in a similar way the Mackey ranges $(\Bi_+, \nui_+) := (\Omegai_+, \bbPi_+)\sslash \Si_+$ and $(\Bi_-, \nui_-) := (\Omegai_-, \bbPi_-)\sslash \Si_-$ associated to the backward cocycle $\chi$. 
	\begin{lem}\label{lem isom Poisson bd chi phi}
		The Mackey ranges $(\Bi_+, \nui_+)$ and $(B_-, \nu_-)$ are isomorphic. The same applies to the Mackey ranges $(\Bi_-, \nui_-)$ and $(B_+, \nu_+)$. 
	\end{lem}
	\begin{proof}
		This is an immediate consequence of Lemma~\ref{lem isom traj chi -}. 
	\end{proof}
	We have the following properties:
	\begin{prop}
		\begin{enumerate}
			\item If the $T$-action is ergodic, then the $G$-action on $B$ is ergodic \cite[Proposition 4.2.24]{zimmer84} (originally due to Mackey \cite{mackey62}). 
			\item The Mackey range $(B, \nu)$ is a cohomology invariant of the cocycle $\varphi$ \cite{zimmer78}. 
		\end{enumerate}
	\end{prop}
	
	The following crucial definition is given in  \cite{bader_furman14}. Recall that for $\eta$ a measure, we denote by $[\eta]$ its measure class. 
	
	\begin{Def}
		We say that $\bnd_+ : \tilde{\Omega} \rightarrow B_+$ and $\bnd_- : \tilde{\Omega} \rightarrow B_- $ are \emph{weakly independent}, or satisfy the \emph{asymptotic past and future independence condition} (apafic), and write $B_-\perp~B_+$ if $(\bnd_- \times \bnd_+)_\ast [\tilde{\mathbb{P}} \otimes \tilde{\mathbb{P}}] $ and ${\bnd_-}_\ast [\tilde{\mathbb{P}}] \otimes {\bnd_+}_\ast [\tilde{\mathbb{P}}]$ represent the same measure class. In other words, the measurable map 
		\begin{align}
			\pi_- \times \pi_+ : (B, \nu) \rightarrow (B_- \times B_+, \nu_- \otimes \nu_+) \nonumber
		\end{align}
		is measure-class preserving. 
	\end{Def}

	The next proposition shows that the definition of apafi ensures that the RDS is admissible in the group, i.e. is not supported on a proper closed subgroup. In terms of Markov chains, this shows ``irreducibility'' of the cocycle. 
	
	\begin{prop}\label{prop apafi admissible}
		Let $(\Omega, \mathbb{P})$ a standard measurable space, $T : \Omega \rightarrow \Omega $ be an ergodic, pmp transformation and $f : \Omega \rightarrow G$ be a measurable map to a discrete countable group $G$. Assume that the cocycle $\varphi : \mathbb{Z} \times \Omega \rightarrow G$ given by $f$ satisfies the asymptotic past and future independence condition. Then $\supp(f_\ast \mathbb{P})$ generates $G$ as a semigroup. 
	\end{prop}
	
	\begin{proof}
		Assume that the semigroup generated by $\supp(f_\ast \mathbb{P})$ is a proper closed subgroup $G_0 < G$. Recall that $\tilde{\Omega}$ describes the stochastic walks of the form $(\varphi (n, \omega)g^{-1})_n$, for $g \in G$. Consider the space $G_0\backslash G$ with the quotient topology. In this case, the trajectory $(\varphi(n, \omega))_n$ is supported on $G_0$, so there is a natural $G$-equivariant, measurable quotient $\tilde{\Omega} \rightarrow G_0\backslash G$. This map passes to the quotients $\tilde{\Omega} \sslash \mathbb{Z} = B$, $\tilde{\Omega}\sslash \mathbb{Z}_{\geq 0} = B_{+}$ and $\tilde{\Omega}\sslash \mathbb{Z}_{\leq 0}= B_{-}$. Hence $G_0\backslash G$ is a common $G$-equivariant quotient of $B_+, B_-$ and $B$, given by measurable maps $\proj_- : B_- \rightarrow G_0\backslash G$, $\proj_+ : B_+ \rightarrow G_0\backslash G$ and $\proj : B \rightarrow G_0\backslash G$. However, this is impossible because $B_+ \perp B_-$. Indeed, take any $F_1, F_2 \subseteq G_0\backslash G$ measurable such that $F_1 \cap F_2 = \emptyset$, and $\nu_-(\proj_{-}^{-1}(F_1)) > 0$, $\nu_{+}(\proj_{+}^{-1}(F_2)) > 0$. Write $A_- = \proj_{-}^{-1}(F_1)$ and $A_+ = \proj_{+}^{-1}(F_2)$. Then $\nu ((\pi_{-} \times \pi_+)^{-1}(A_- \times A_+)) = 0$, a contradiction. 
	\end{proof}
	
	\subsection{Equivariant maps for ideal pasts and futures}
	Let $(\Omega, \mathcal{F}, \mathbb{P}, T, f, G)$ be a discrete RDS on a locally compact second countable group $G$ with $T$ invertible and measure-preserving.

	Let now $V$ be a Polish space on which $G$ acts by Borel transformations. 
	It is enlightening to consider the case where $V = \prob(X)$, for $X$ a Polish $G$-space. In the case where $\varphi$ is the random walk $(Z_n(\omega)= \omega_1\dots \omega_{n})_n$ on $G$ associated to an admissible measure $\mu$, a standard strategy to study the right random walk $(Z_n(\omega )o)_n$ starting at a given point $o \in X$ is to find $\mu$-stationary measures on $X$, that is, probability measures $\nu$ that satisfy $\mu \ast \nu = \nu $. One of the fundamental results of boundary theory relates stationary measures, $T$-equivariant measurable maps $\Omega \to \prob(X)$ and $G$-equivariant measurable maps from the Poisson-Furstenberg boundary $B=B(G, \mu)$ to $\prob(X)$. 
	
	\begin{thm}[{\cite{furstenberg63}}]\label{thm dual furst}
		Let $(B(G, \mu), \nu_B)$ be the Poisson-Furstenberg boundary associated to the probability measure $\mu$ on $G$, and let $X$ be a Lebesgue $G$-space. Then the following assertions are equivalent: 
		\begin{enumerate}
			\item the probability measure $\nu \in \prob(X)$ is $\mu$-stationary;
			\item \label{equivariance map classical}there exists a measurable map $\Phi : \Omega \to \prob(X)$ such that $\mathbb{P}$-almost surely, $\Phi(\omega) = (\omega_1)_\ast \Phi(T\omega)$ and 
			$$\nu = \int \Phi (\omega) d \mathbb{P}(\omega)$$
			\item there exists a measurable $G$-equivariant map $ \phi : B \rightarrow \prob(X)$ such that $$\nu = \int \phi (b) d \nu_B(b).$$
		\end{enumerate}
	\end{thm}
	
	As we have seen, the analogue of stationary measures in the broader context of random dynamical systems is the set of $(T^+)^X$-forward invariant measures. The following theorem is the RDS extension of Theorem \ref{thm dual furst}. Let $V $ be a  compact convex metrizable space $V$ on which $G$ acts by affine transformations. To make the comparison as explicit as possible, we define the space of measurable maps
	\[\Map_f(\Omegai_+, V):=  \{ \phi  \in \Map(\Omega_-, V) \mid  \phi(\omega)= f(\Ti\omega)\phi(\Ti \omega) \text{ almost everywhere}\} \nonumber.\]

	\begin{thm}\label{thm corresp t_inv G_inv}
		Let $(\Omega, \mathcal{F}, \mathbb{P}, T, f, G)$ be a an invertible discrete RDS on a locally compact second countable group $G$, and let $V$ be a compact convex metrizable space on which $G$ acts by affine transformations (e.g. $V= \prob(X)$ for $G\curvearrowright X$ a continuous action by isometries on a compact metrizable space $X$). There are natural bijections 
		$$\Map_G(\Bi_+, V) \simeq \Map_f(\Omegai_+, V), \hspace{2em} \Map_G(\Bi_-, V) \simeq \Map_f(\Omegai_-, V),$$ 
		$$ \Map_G(B_-, V) \simeq \Map_f(\Omega_-, V), \hspace{2em} \Map_G(B_+, V) \simeq \Map_f(\Omega_+, V)$$ 
		and
		$$\Map_G(B, V) \simeq \Map_f(\Omega, V) \simeq \Mapstat(\Omega_+, V).$$
		
	\end{thm}
	
	\begin{proof}
		The bijection $\Map_f(\Omega, V) \simeq \Mapstat(\Omega_+, V)$ is the bijection given by Theorem \ref{thm lim measures rds}. Let us prove that $\Map_G(\Bi_+, V) \simeq \Map_f(\Omegai_+, V)$. As $(\Bi_+, \nui_+) = (\check{\Omega}_+, \bbP_+)\sslash \Si_+$, the space $\Map_G(\Bi_+, V) $ identifies with the space $\Map_G(G^{\bbN}, V)^{\Si_+}$ of $\Si_+$-invariant, $G$-equivariant maps from $G^{\bbN}$ to $V$. Let $\Psi \in \Map_G(\Bi_+, V)$ be such a map. 
		Recall that by Equation~\eqref{eq intertwining chi}, for all $m \geq 0$ and $\omega \in \Omegai_+$, 
		\begin{align}
			\Si_+^m\alphai^+_g(\omega) = \alphai^+_{g\chi(m, \omega)} (\Ti_+^m \omega). \label{eq intertwining chi +}
		\end{align}
		Then we have
		\begin{align}
			\Psi (\alphai^+_g (\omega)) &= \Psi (S_+\alphai^+_g (\omega)) \text{ by $\Si_+$-invariance} \nonumber\\
			&= \Psi (\alphai^+_{gf(\Ti \omega)} (\Ti\omega)) \text{ by Equation \eqref{eq intertwining chi +}} \nonumber\\
			& =  \Psi (gf(\Ti \omega) g^{-1} \alphai^+_g (\Ti\omega)) \nonumber\\
			& =  gf(\Ti \omega) g^{-1} \Psi (\alphai^+_g (\Ti\omega)) \text{ by $G$-equivariance}. \label{eq equiv chi +}
		\end{align}
		Using Fubini's Theorem, we let $\psi : \Omegai_+ \to V $ by $\psi (\omega ) = g^{-1}\Psi(\alphai_g^+(\omega))$ for almost every $g$. It follows from the construction that $\psi$ is essentially well-defined and does not depend upon the choice of the representative $\Psi$. By Equation~\eqref{eq equiv chi +}, we obtain that 
		\[\psi(\omega) = f(\Ti \omega ) \psi(\Ti \omega). \]
		Therefore, we get a well-defined map 
		$$\Map_G(\Bi_+, V) \to \Map_f(\Omegai_+, V).$$ 
		
		Now let $\psi \in \Map_f(\Omegai_+, V)$ and denote again by $\psi$ a representative. By construction, $\alphai_+$ defines an isomorphism between $(\Omega_- \times G, \mathbb{P}^-\otimes \mu_G) $ and $(\Omegai_+ = G^{\bbN}, \bbPi_+ )$. Consequently, there exist essentially well-defined measurable map $a : \Omegai_+ \to \Omega_- $ and $ b : \Omegai_+ \to G$ such that for almost every $\omegai= (g_0, g_1, \dots)$, 
		$$\alphai^+_{b(\omegai)} (a(\omegai)) = \omegai.$$
		Observe that by construction, the map $b$ is defined by $b(g_0, g_1, \dots) = g_0$ and is then $G$-equivariant. 
		Define $\Psi: G^{\bbN} \to V$ by $\Psi (\omegai)= b(\omegai)\psi (a(\omegai))$. Then $\Psi (g \omegai)= gb(\omegai)\psi (a(\omegai))$ and $\Psi $ is almost everywhere $G$-equivariant. Remark that by construction, for almost every $\omegai= (g_0, g_1, \dots)$, 
		$$ g_1 = b(\Si \omegai) = (\alpha^+_{g_0}(a(\omegai)))_1 = g_0 f(\Ti(a(\omegai))).$$
		Observe also that $a(\omegai)$ is defined in order to verify that $\alphai^+_e(a(\omegai)) = \omegai = (\chi(n, a(\omegai)))_n$, so that we have the identity: for $g_0 \in G$, $\omega \in \Omega_-$, 
		\begin{align}
			\Psi(\alphai^+_{g_0}(\omega)) = g_0 \psi(\omega) \label{eq identity a}.
		\end{align}
		Using again Equation \ref{eq intertwining chi}, we get that for almost every $\omegai= (g_0, g_1, \dots)$,
		\begin{align}
			\Psi(\Si_+ \omegai) & =  \Psi (\Si_+ \alpha^+_{b(\omegai)} (a(\omegai))) \nonumber \\
			& =  \Psi (\alpha^+_{b(\omegai) f(\Ti a(\omegai))} (\Ti a(\omegai))) \nonumber \\
			& =  b(\omegai) f(\Ti a(\omegai)) \Psi(\alphai^+_e(\Ti a(\omegai))) \nonumber \\
			& = b(\omegai) f(\Ti a(\omegai))\psi(\Ti a(\omegai)) \text{ by Equation~\eqref{eq identity a}}\nonumber\\
			& =  b(\omegai) f(\Ti a(\omegai)) f(\Ti a(\omegai))^{-1}\psi(a(\omegai)) \text{ because $\psi \in \Map_f(\Omegai_+, V)$} \nonumber\\
			& =  \Psi( \omegai) \nonumber.			
		\end{align}
		Therefore $\Psi$ is $G$-equivariant and $\Si_+$-invariant, hence belongs to $\Map_G(\Bi_+, V)$. This procedure defines an inverse to the map $$\Map_G(\Bi_+, V) \to \Map_f(\Omegai_+, V)$$ defined above. 
		
		The proofs that $\Map_G(B_+, V) \simeq \Map_f(\Omega_+, V)$ and $\Map_G(B_-, V) \simeq \Map_f(\Omega_-, V)$ are completely analogous. 
	\end{proof}

	Consider the inclusion $\Map_f(\Omega_-, V) \subseteq \Map_f(\Omega, V)$ given by the natural projection $\Omega \to \Omega_-$. By the previous theorem, these maps correspond to the non-empty compact convex set of maps $\Map_G(B_-, V) $, and thus form a non-empty convex subset of $\Mapstat(\Omega_+, V)$. 
	\begin{Def}
		We call the subset of $\Mapstat(\Omega_+, V)$ corresponding to the maps $\phi \in \Map_f(\Omega_-, V)$ the set of \emph{past-oriented stationary maps}, and denote them by $\Mapstat_-(\Omega_+, V)$. 
	\end{Def}
	The previous discussion can be summarized in saying that the following diagram commutes with respect to the natural inclusions and isomorphisms:
	
	\[
	\begin{tikzcd}
		\Map_G(B_-, V) \arrow[hookrightarrow]{d}{\pi_{-}^{\ast}} \arrow[r,  "\sim"] & \Map_f(\Omega_-, V) \arrow[hookrightarrow]{d}{p_-^\ast} \arrow[r, " \sim "]& \Map_{-}^{\text{stat}}(\Omega_+, V) \arrow[hookrightarrow]{d} \\
		\Map_G(B, V) \arrow[r,  "\sim "]  & \Map_f(\Omega, V) \arrow[r,  "\sim "] & \Mapstat(\Omega_+, V)  \\
		\Map_G(B_+, V) \arrow[hookrightarrow]{u}{\pi_{+}^{\ast}} \arrow[r,  "\sim "] & \Map_f(\Omega_+, V) \arrow[hookrightarrow]{u}{p_{+}^{\ast}} &
	\end{tikzcd}
	\]
	The vertical arrows follow from the constructions. If $V$ is compact, $\Map_G(B_-, V)$ is non-empty and all these sets are then non-empty.

	\subsection{Boundary systems and ergodic properties}

	In this section, we introduce the notion of $G$-boundaries and boundary systems. To summarize, $G$-boundaries are measured spaces that satisfy the ergodic properties of the Poisson-Furstenberg boundary associated to an admissible measure. We refer to \cite{bader_furman14}, where the vocabulary of fiberwise isometric actions was introduced. 
	
	Throughout this section, $G$ is a locally compact second countable group. Recall that a standard probability space $(S, \eta)$ is a $G$-space if there is a non-singular action $G \curvearrowright S$, and a map $\psi : (S, \eta) \to (S', \eta')$ between $G$-spaces is called a $G$-map if it is a $G$-equivariant measurable map. 
	
	\begin{Def}
		Let $(S, \eta)$ be a $G$-space where $\eta \in \prob(S)$. We say that the $G$-action on $S$ is \emph{isometrically ergodic} if for all separable metric space $(Y, d)$ on which $G$-acts by isometries, and for all $G$-equivariant measurable map $\psi: S \rightarrow Y$, alors $\psi$ is essentially constant. In other words, for any such $(Y,d)$, there is no non-trivial $G$-map  
		\[
		\begin{tikzcd}
			\psi : (S, \eta ) \arrow[r, "G"]  & (Y, d). 
		\end{tikzcd}
		\]
		If the diagonal action of $G$ on $(S \times S, \eta \times \eta)$ is isometrically ergodic, we say that $G \curvearrowright S$ is \emph{doubly isometrically ergodic}. 
	\end{Def}
	
	In other words, the $G$-action on $(S, \eta)$ is isometrically ergodic if for all metric space $(Y, d)$ on which $G$-acts by isometries,
	\begin{align}
		\Map_G(S, Y) = \Map_G(\{\ast\}, Y) \nonumber
	\end{align}
	where we put the Borel $\sigma$-algebra on $(Y,d)$ given by the metric. 
	
	Let us now define the notion of coarse metric ergodicity, introduced in \cite{bader_caprace_furman_sisto22}.
	
	If  $(Y, \mathcal{Y}, \nu)$ is a standard probability space, an \emph{essentially well-defined coarse metric} on $Y$ is the datum of a $\nu$-conull subset $Y' $ and a coarse metric $d: Y' \times Y' \to \bbR_+$. 
	\begin{Def}
		Let $G$ be a locally compact second countable group. We say that the $G$-action on a standard probability space $(Y, \mathcal{Y}, \nu)$ is \emph{coarsely metrically ergodic} if every essentially well-defined and $G$-invariant Borel coarse metric $d$ is essentially bounded, that is, there exists $R > 0$ such for a.e. pair $(x, y) \in \Omega\times \Omega$, $d(x, y ) \leq R$ if $d(x, y )$ is defined. 
	\end{Def}
	\begin{rem}
		Equivalently, the $G$-action on $\Omega$ is coarsely metrically ergodic if every $G$-invariant Borel metric on any conull subset $\Omega_0 \subseteq \Omega$ as above is essentially bounded, see \cite[Remarck 2.2]{bader_caprace_furman_sisto22}. One can also observe that if $G$ is locally compact second countable, any $G$-invariant coarse metric on an ergodic $G$-space $\Omega$ is coarsely separable, i.e. $\Omega$ can be covered by a countable collection of balls of a fixed radius, see \cite[Remark 2.6]{bader_caprace_furman_sisto22}. 
	\end{rem}

	We say that $(A, \eta)$ is an \emph{amenable} $G$-space if $(A, \eta)$ is a standard Borel space such that there exists a non-singular and amenable action $G \curvearrowright (A, \eta)$ in Zimmer's sense, see \cite{zimmer84} and \cite[Section 7.1]{duchesne13}.

	\begin{Def}
		Let $G$ be a locally compact second countable group and let $(B, \nu)$, $(B_-, \nu_-)$ and $ (B_+, \nu_+)$ be $G$-Lebesgue probability spaces with $G$-equivariant measurable maps $\pi_+ : B \rightarrow B_+ $ and $ \pi_- : B \rightarrow B_-$. We say that $(B, B_-, B_+)$ forms a \emph{$G$-boundary system} if 
		\begin{enumerate}
			\item $\pi_- \times \pi_+ : (B, \nu) \rightarrow (B_- \times B_+, \nu_- \otimes \nu_+)$ is measure-class preserving; 
			\item The $G$-actions on $B$, $B_-$, and $B_+$ are Zimmer-amenable;
			\item $\pi_+ : B \rightarrow B_+ $ and $ \pi_- : B \rightarrow B_-$ are relatively isometrically ergodic (in the sense of \cite{bader_furman14}). 
		\end{enumerate}
	\end{Def}

\subsection{Ergodic properties of the Mackey range of cocycles}\label{section erg prop mackey range}
	
	Coming back to the RDS setting $(\Omega, \mathcal{F}, \mathbb{P}, T, f, G)$, and taking the space of ergodic components (the Mackey range of the cocycle $\varphi$) as defined above, we obtain standard probability $G$-spaces $(B, \nu)$, $(B_+, \nu_+)$ and $(B_-, \nu_-)$. If $T$ is invertible ergodic, we saw that the $G$-action on any of these spaces is ergodic. In fact we can say much more.

	\begin{thm}[{Coarse isometric ergodicity of the Mackey range}]\label{thm coarse isom erg}
		Let $G$ be a locally compact second countable group and let $(\Omega, \mathcal{F}, \mathbb{P}, T, f, G)$ be an invertible RDS  on $G$ satisfying the apafi condition. Let $\Bi_-$ and $\Bi_+$ be the corresponding past and future factors associated to $\chi$. Then the diagonal $G$-action on $\Bi_- \times \Bi_+$ is metrically ergodic and coarsely metrically ergodic. 
	\end{thm}

	The goal of this section is to prove this result. Consider the \emph{Poisson transform} associated to $\varphi$:
	\[F \in L^\infty(\Bi_+, [\nui_+]) \mapsto \Phi_\chi(F) \in L^\infty(G, [\mu_G]),\]
	where 
	\[\Phi_\chi(F) (g) = \int_{\Bi_+} F(gb)d\nui_+ (b)\]
	for all $g \in G$. Notice that $\Phi_\chi$ is a linear positive contraction, unital once we fix $\nui_+$ to be a probability measure. 
	\begin{prop}\label{prop Poisson transform}
		Let $F \in L^\infty(\Bi_+, [\nui_+])$. Then almost surely, the limit 
		\[\hat{F}(\omega):= \underset{n \to \infty}{\rightarrow} \Phi_\chi(F) (\chi(n, \omega))\]
		exists and the function $\hat{F} $ is in $L^\infty( \Omega_-, [\bbP_-])$. If we identify $L^\infty( \Omega_-, \bbP_-)$ with $L^\infty(\Omegai_+, [\bbPi_+])$, $\hat{F} $ is $\Si_+$-invariant so that $\hat{F}\in L^\infty(\Bi_+, [\nu_+])$ and $\hat{F} = F $ almost surely. 
	\end{prop}
	\begin{proof}
		We denote by $\hat{F}_n (\omega) = \Phi_\chi(F) (\chi(n, \omega))$, and consider the increasing sequence of $\sigma$-subalgebras $\calF_{-n}^{0}$ from Section \ref{section past future sigma algebras}, so that $\calF_{-\infty}^{0}= \calF^-$. For all $n $, $\|\hat{F}_n\|_\infty \leq \|F\|_\infty$. Then 
		\begin{align}
			\mathbb{E}(\hat{F}_{n+1} \mid \calF_{-n}^{0}) (\omega) &= \int_G \Phi_\chi(F) (\chi(n, \omega) g)d(f\circ T^{-(n+1)})_\ast \bbP(g) \nonumber \\
			&= \int_G \int_{B_+} F(\chi(n, \omega) gb)d(f\circ T^{-(n+1)})_\ast \bbP(g)d\nu_+(b) \nonumber \\
			&= \int_\Omega \int_{B_+} F(\chi(n, \omega) f(T^{-(n+1)}\omega')b)d\bbP(\omega')d\nu_+(b) \nonumber \\
			&=  \int_{B_+} F(\chi(n, \omega) b')d\nu_+(b'), \nonumber 
 		\end{align}
		where we used the fact that $F \in L^\infty(B_+, [\nu_+]) $ can be viewed as a $\Si_+$ invariant function in $ L^\infty (\Omegai_+, \bbPi_+)$ to pass from the third line to the fourth. Therefore $\hat{F}_n$ is a uniformly bounded martingale, hence it converges $\bbP_-$-almost everywhere to an essentially well-defined $\hat{F}\in L^\infty( \Omega_-, \bbP_-)$, with $\|\hat{F}\|_\infty \leq \|F\|_\infty$. 
		
		Identify $L^\infty( \Omega_-, \bbP_-)$ with $L^\infty (\Omegai_+, \bbPi_+)$ via $\alphai_e^+$ as in Section \ref{section space trajectories}, so that $\hat{F}$ is seen as a map in $L^\infty(\Omegai_+, \bbPi_+)$. 
		By definition, $\hat{F}$ is $\Si_+$-invariant, therefore  $\hat{F} \in L^\infty(B_+, \nu_+)$.
		
		Finally, observe that 
		\begin{align}
			\mathbb{E}(\hat{F} \mid \calF_{-n}^0) (\omega) &= \Phi_\chi(F)(\chi(n, \omega)) \nonumber .
		\end{align}
		As the filtration $\calF_{-n}^0$ satisfies that $\calF_{-\infty}^0 = \calF^-$, we obtain that $\hat{F} = F$ almost everywhere, finishing the proof. 
	\end{proof}
	
	The strategy for the proof of Theorem \ref{thm coarse isom erg} relies on the following lemma. For a measurable subset $E \subseteq A \times A'$  in a product of measurable spaces, and for $a \in A$, we define 
	$$E_a := \{ a' \in A' \mid (a,a') \in E\} $$
	the right coset in $E$ determined by $a$. It is a measurable subset of $A'$. 
	
	\begin{lem}\label{lem sat rec}
		Let $E \subseteq \Bi_- \times \Bi_+$ be a measurable subset of positive $\nui_- \times \nui_+$-measure, and let $\varepsilon >0$. Denote by $\pi_\pm : \Bi_-\times \Bi_+ \to \Bi_\pm$ the natural projection on the corresponding factor. Then for almost every $b_- \in \Bi_-$, there exists $g \in G$ such that 
		\begin{enumerate}
			\item (SAT) $\nui_+(gE_{b_{-}}) >1- \varepsilon$ ; 
			\item (Recurrence) $gb_- \in E_-:= \pi_-(\Bi_- \times \Bi_+)$. 
		\end{enumerate}
	\end{lem}
	This lemma roughly states that the $G$-action on $B_+$ is SAT, and that one can simultaneously use Poincar\'e recurrence theorem on the other factor. 
	\begin{proof}[Proof of Lemma \ref{lem sat rec}]
		Recall that $\Bi_- $ and $\Bi_+$ come from the factors 
		$$ \check{\bnd}_+ : \Omegai_+ \rightarrow \Bi_+ \hspace{2em} \text{ and }\hspace{2em} \check{\bnd}_- : \Omegai_- \rightarrow \Bi_-. $$ 
		Using the projections $\Omegai \to \Omegai_\pm $ and the measurable isomorphism $\alphai_e : \Omega \to \Omegai$ defined in Section \ref{section space trajectories}, We can consider the preimage $\tilde{E} $ in $\Bi_- \times \Omega$. In order to ease the notations, we also denote by $\pi_- : \Bi_- \times \Omega  \to \Bi_-$ the projection on the first factor.  Define the transformation $\mathcal{S} $ on $\Bi_- \times \Omega_-$ by 
		$$\mathcal{S}(b_-,\omega) = (f(\Ti\omega)^{-1} b_-, \Ti\omega),$$
		so that $$\mathcal{S}^n (b_-,\omega) = (\chi(n, \omega)^{-1}b_-, \Ti^n \omega).$$ By construction of $\Bi_-$, the transformation $\tilde{S}$ is measure-preserving. By Poincar\'e recurrence theorem, for almost every $\tilde{e} \in \tilde{E}$, there are infinitely many $n \in \bbN$ such that $\mathcal{S}^n \tilde{e}\in \tilde{E}$. 
		
		Notice that by Fubini, for $\nui_-$-almost every $b_- \in \pi_-(E)$, $\nui_+(E_{b_{-}}) >0 $. Equivalently, for almost every $\tilde{e}\in \tilde{E}$, $\bbP_-(E_{\pi_{-}(\tilde{e})}) >0 $.

		Consider $b_- = \pi_- (\tilde{e}) $ for $\tilde{e}$ in the intersection of the two conull sets above, so that $E_{b_-}$ is of positive $\nui_+$-measure. Let $\Phi_\chi(\mathds{1}_{{E_{b_-}}})$ be its Poisson transform defined by 
		$$\Phi_\chi(\mathds{1}_{{E_{b_-}}}) (g) = \int_{\Bi_+} \mathds{1}_{E_{b_-}} dg_\ast \nui_+(x). $$
		Then by Proposition \ref{prop Poisson transform}, for every $\omega \in \tilde{E}_{b_-}$, one has the convergence: 
		
		$$\Phi_\chi(\mathds{1}_{{E_{b_-}}}) (\chi(n, \omega)) \to \mathds{1}_{{E_{b_-}}}(\omega) = 1.$$
		In other words, for every $\omega \in \tilde{E}_{b_-}$, 
		\begin{align}
			\nui_+ ( (\chi(n, \omega))^{-1}E_{b_-} ) \underset{n \to \infty}{\to} 1. \label{eq SAT}
		\end{align}
		
		In particular, for any $\varepsilon > 0 $, and for almost every $b_-$, there exists a well-chosen $\omega \in \Omega $ and $n$ large enough so that $\nui_+ ( (\chi(n, \omega))^{-1}E_{b_-} ) > 1 - \varepsilon$ and $\chi(n, \omega)^{-1} b_- \in E_-$, which finishes the proof. 
	\end{proof}
	\begin{proof}[Proof of Theorem \ref{thm coarse isom erg}]
		The rest of the proof is now exactly the same as in \cite[Proof of Theorem 2.7]{bader_caprace_furman_sisto22}, which only relies on Lemma \ref{lem sat rec} above and general measure-theoretic manipulations.
	\end{proof}
	
	The following definition was introduced by Jaworski \cite{jaworski95}. 
	\begin{Def}
		Let $\nu'$ be a $\sigma$-finite measure on a Borel space $G$-space $Y$, and assume that $\nu'$ is $G$-quasi-invariant. We say that the $G$-space $(Y, \nu)$ is \emph{strongly approximately transitive} if there is a probability measure $\nu \ll \nu' $ absolutely continuous with respect to $\nu'$ such that the convex hull of the orbit $G\nu$ in the space $L^1(Y,\nu')$ is dense in the space of absolutely continuous probability measures on $Y$ for the topology induced by the total variation norm. 
	\end{Def}

	\begin{cor}\label{cor sat}
		The $G$-action on $\Bi_\pm$ is SAT. 
	\end{cor}
	\begin{proof}
		By \cite[Proposition 2.2]{jaworski95}, the $G$-action on $(\Bi_\pm, \nui_\pm)$ is SAT if and only if it admits a probability $\nui_\pm' \ll \nui_\pm$ such that for every $A $  measurable subset of $\Bi_\pm$ of positive $\nui_\pm$-measure, and every $\varepsilon >0$, there exists $g \in G$ such that $\nui'_+(gA) \geq 1- \varepsilon$. This fact was shown in the proof of Lemma \ref{lem sat rec}, in Equation \eqref{eq SAT}. 
	\end{proof}

	Let us give another ergodic property of the spaces $\Bi, \Bi_\pm$, which is stated in \cite[Theorem 4.1]{bader_furman25}.

	\begin{thm}\label{thm boundary system}
		Let $(\Omega, \mathcal{F}, \mathbb{P}, T, f, G)$ be a greg on a locally compact second countable group satisfying the past and future asymptotic condition. Assume that the flow $T$ is ergodic. Then $(\Bi, \nui)$, $(\Bi_+, \nui_+)$ and $(\Bi_-, \nui_-)$ form a boundary system. 
	\end{thm}
	
	\begin{proof}
		By the apafi condition, $\pi_- \times \pi_+ : (B, \nu) \rightarrow (B_- \times B_+, \nu_- \otimes \nu_+)$ is measure-class preserving, and so the same statement holds for $\Bi, \Bi_\pm$ by Lemma~\ref{lem isom Poisson bd chi phi}. The theorem is then a combination of Lemma \ref{lem Mackey range amenable} below and \cite[Theorem~3.3]{bader_furman25}. 
	\end{proof}

	\begin{lem}\label{lem Mackey range amenable}
		The $G$-spaces $\Bi, \Bi_+, \Bi_-$ are Zimmer-amenable. 
	\end{lem}
	\begin{proof}
		We only treat the case for $\Bi_+$, the case $\Bi_-$ being completely analogous. We have to show that for every affine action on a non-empty metrizable compact convex space $Q$, $\Map_G(\Bi_+, Q) \neq \emptyset$. By Theorem \ref{thm corresp t_inv G_inv}, this is equivalent to showing that $\Map_f(\Omegai_+, Q)\neq \emptyset$. As the set $Q$ is non-empty, $\Map(\Omegai_+, Q)\neq \emptyset$ and we can apply Theorem \ref{thm krylov}. 
		
		Therefore $\Map_G(\Bi_+, Q)$ and $\Map_G(\Bi_-, Q)$ are non-empty, so $\Map_G(\Bi, Q)$ is non-empty as well. 
	\end{proof}

	As an immediate consequence of Lemma~\ref{lem isom Poisson bd chi phi}, we obtain analogous statements for the $G$-spaces $B, B_-, B_+$. 
	\begin{cor}
		Theorems~\ref{thm coarse isom erg} and \ref{thm boundary system}, Lemma~\ref{lem Mackey range amenable} and Corollary~\ref{cor sat} remain true if one replaces $\Bi, \Bi_+, \Bi_-$ by $B,B_-,B_+$ respectively. 
	\end{cor}

	\subsection{Atomic parts of stationary maps}
	We keep the notations introduced at the beginning of the section. In particular, we let $(\Omega, \mathcal{F}, \mathbb{P}, T, f, G)$ be a greg satisfying the apafi condition on a locally compact second countable group $G$, and we assume that the invertible measure-preserving map $T$ is ergodic. In particular, by Theorem \ref{thm boundary system}, the spaces of ideals $(B, B_-, B_+)$ and their natural probability measure classes form a boundary system for $G$. We let $(X, \mathcal{B})$ be a compact metrizable space on which $G$ acts by homeomorphisms, with $\calB$ its $\sigma $-algebra of Borel sets. 
	
	Recall that $\prob(X)$ is endowed with the topology of weak-$\ast$~convergence, making it a convex compact metrizable space on which $G$ acts by affine transformations. The following result shows that past-oriented stationary measures are diffuse. 
	
	\begin{prop}[{\cite[Proposition 4.7]{bader_furman25}}]\label{prop purely non atomic}
		If the action of $G$ on $X$ does not have finite orbits, then for any $\mu^+_\cdot \in \Mapstat_-(\Omega_+, \prob(X))$, the probability measure $\phi_+(b_+)$ is almost surely purely non-atomic. 
	\end{prop}
	
	\section{Ergodic cocycles in Gromov-hyperbolic spaces}\label{section erg coc hyp cv}

	\subsection{Boundary maps in hyperbolic spaces}\label{section bd map hyp}

	In this section, let $G$ be a discrete countable group. Let $(X,d)$ be a separable almost-geodesic Gromov-hyperbolic space and $G \curvearrowright (X, d)$ be a non-elementary action by isometries. Recall that $\Xh$ denotes the horofunction compactification of $X$.

	As we assume $X$ separable, $\Xh$ is a compact metrizable space on which $G$ acts by homeomorphisms. Therefore, by amenability of the actions $G \curvearrowright B, B_\pm $, the set $\Map_G ( B_-, \prob(\Xh)) $ is non-empty. 
	
	The following result is stated for geodesic hyperbolic spaces in \cite{bader_caprace_furman_sisto22} but, with the preparations made in Section \ref{section prelim hyp} and Lemma \ref{lem coarse metric hyp}, it remains valid for almost geodesic hyperbolic spaces. 
	
	\begin{thm}[{\cite[Theorem 3.1]{bader_caprace_furman_sisto22}}]\label{thm uniq bdry map hyp}
		Let $(B, B_-, B_+) $ be any $G$-boundary system, and $G \curvearrowright X$ as above. Then there exists $G$-maps $\phi_- : B_- \rightarrow \bdg X$ and $\phi_+ : B_+ \rightarrow \bdg X$ such that the map 
		\begin{align}
			\phi_{\bowtie}: (b_-, b_+) \in B_-\times B_+ \mapsto (\bdg X)^2 \nonumber
		\end{align} 
		is essentially contained is the set of distinct pairs of points of the boundary $\bdg(X)^{(2)}$. Moreover: 
		\begin{enumerate}
			\item[(i)] $\Map_G ( B_-, \prob(\bdg X)) = \{\delta \circ \phi_-\}$ and $\Map_G ( B_+, \prob(\bdg X)) = \{\delta \circ \phi_+\}$. 
			\item[(ii)] $\Map_G ( B_- \times B_+, \bdg X) = \{\phi_-\circ \text{pr}_-, \phi_+\circ \text{pr}_+\}$ 
			\item[(iii)] $\Map_G ( B_- \times B_+, \bdg X^{(2)}) = \{\phi, \tau\circ\phi\}$, where $\tau (\xi,\eta) = (\eta, \xi)$. 
		\end{enumerate} 
	\end{thm}
	
	This description will be applied to the boundary system $(B, B_\pm)$ induced by a discrete RDS on $X$ associated to a $G$-action. 
	
	\subsection{Convergence of the ergodic cocycle}\label{section rw hyp}
	For the rest of this section, $(\Omega, \mathcal{F}, \mathbb{P}, T, f, G)$ is a greg on a discrete countable group $G$, and assume that the measure-preserving map $T$ is invertible and ergodic. We assume that the RDS $\varphi : \mathbb{Z} \times \Omega \rightarrow G$ given by $f$ satisfies the asymptotic past and future independence condition. In particular, by Theorem \ref{thm boundary system}, the spaces of ideals $(B, B_-, B_+)$ and their natural probability measure classes form a boundary system for $G$. 
	Let $(X,d)$ be a separable almost-geodesic Gromov-hyperbolic space and $G \curvearrowright (X, d)$ be a non-elementary action by isometries.

	As in Section \ref{section rds}, we denote by $\chi$ the backward random cocycle defined by 
	$$\chi(n, \omega) = \varphi(-n, \omega)^{-1} = \varphi(n, T^{-n}\omega)$$
	for $n \in \bbZ$, $\omega \in \Omega$. 
	
	The main result of this section is the following convergence theorem. 
	
	\begin{thm}\label{thm cv cocycle hyp}
		For any $o\in X$, the backward cocycle $(\chi(n, \omega) o)_n$ converges $\mathbb{P}$-almost surely to a point $\xi^+(\omega) \in \partial X$ in the Gromov boundary. 
	\end{thm}
	
	The next proposition combines Theorem~\ref{thm uniq bdry map hyp} with Theorem~\ref{thm lim measures rds}.

	\begin{prop}\label{prop lim measures hyp}
		There exists a measurable map $\psi \in \Map_f(\Omega, \bdg X)$ such that for $\bbP$-almost every $\omega \in \Omega$, $\mu_\omega $ is the Dirac measure $\delta_{\psi(\omega)}$ at $\psi (\omega) \in \bdg X$. In particular, we have the almost sure convergence 
		$$ \chi(n, \omega)\mu^+_{p^+(T^{-n}\omega)} \underset{n\to \infty}{\longrightarrow} \delta_{\psi(\omega)}$$
		in the weak-$\ast$ topology. 
	\end{prop}

	\begin{proof}
		As $\Map_G(B_-, \prob(\bdg X))$ is non-empty by Theorem \ref{thm uniq bdry map hyp}, Theorem \ref{thm corresp t_inv G_inv} implies that there exists a past-oriented stationary measure, which we denote by $\mu^+_\cdot$.  Using the one-to-one correspondence of Theorem \ref{thm lim measures rds}, there exists an invariant measure $\mu_\cdot \in \Map_f(\Omega, \prob(\bdg X))$ such that 
		$$\mu_\omega = \lim_n \chi(n, \omega))\mu^+_{p^+(T^{-n}\omega)}$$
		in the weak-$\ast$ topology. 
		Moreover, an immediate combination of Theorems \ref{thm corresp t_inv G_inv} and \ref{thm uniq bdry map hyp} gives that $\mu_\omega$ is almost surely a Dirac measure.
	\end{proof}
	Proposition \ref{prop lim measures hyp} represents a contraction of ``most'' of the boundary to the boundary point $\phi(b)$. 
	The last piece needed for the proof of Theorem \ref{thm cv cocycle hyp} is then a better understanding of the support of the measures $\{\mu^+_\omega\}$. In the random walk case, each one of these measures is the same stationary measure, but here the support $\supp(\mu^+_\omega)$ is a random closed set, as defined in, for instance \cite{crauel02}. 
	\begin{lem}\label{lem control support}
		Let $g \in G$ be a loxodromic isometry, and let $\xi^+, \xi^-$ be its attractive and repelling fixed point. Then $\xi^+$ belongs almost surely to the support of $\mu^+_{p^+ (\omega)}$. 
	\end{lem}
	
	\begin{proof}
		Let $g, \xi^+ $ and $\xi^-$ be as in the statement. Let $U$ be an open neighborhood of $\xi^+$ and let $V$ be an open neighborhood of $\xi^-$. By \cite[Corolloary 1.6.5]{arnold98}, the support $\omega_+ \in \Omega_+  \mapsto \supp \mu^+_{\omega_+} \subseteq \bdg X$ is a random closed set. By \cite[Proposition 1.6.2]{arnold98}, the set 
		$$\Omega(V) := \{\omega \in \Omega \mid \supp\mu^+_{p^+(\omega)} \cap (\bdg X  - \overline{V})\neq \emptyset\} $$
		is a measurable set. By Proposition \ref{prop purely non atomic}, the measure $\mu^+_{p^+ (\omega)}$ is almost surely purely non-atomic. Thus for all $\varepsilon >0$, there exists an open neighborhood $V_\varepsilon$ of $\xi^-$ such that $\bbP (\Omega(V_\varepsilon)) >1- \varepsilon >0 $. Fix such a set. 
		
		By North-South dynamics on the Gromov boundary, there exists $k_0$ such that for all $k \geq k_0$, $g^k(\bdg X  - V_\varepsilon) \subseteq U$. By the asymptotic past and future condition and Proposition \ref{prop apafi admissible}, for almost every $\omega \in \Omega$, there exists $n\in \bbN$ such that 
		$$ \chi(n, \omega ) = g^k,$$ 
		for some $k \geq k_0$. In particular, if we write 
		\[\Omega_n:= \{ \omega \in \Omega\mid \exists k \geq k_0 \colon \chi(n, \omega) = g^k\},\]
		then $\Omega':= \bigcup_n \Omega_n $ is conull in $\Omega$. 
		Let $\omega $ be in this set. 
		
		Let us rewrite the stationarity equation: $ \mu^+_{p^+(\omega)} (U)=0 $ if and only if, for all $n \in \bbN$, 
		\[\int_\Omega \varphi(n, T^{-n}\omega')\mu^+_{p^+(T^{-n}\omega')} (U)d\mathbb{P}^+_\omega(\omega') =0. \]
		In other words,  $ \mu^+_{p^+(\omega)} (U) =0$ if and only if, for a.e. $\omega' \in (p^+)^{-1}(p^+(\omega))$, for all $n \in \bbN$, 
		\[ \chi(n, \omega')\mu^+_{p^+(T^{-n}\omega')} (U)= 0.\]
		In particular, it implies that for all $n$, and for all $\omega' \in \Omega_n \cap (p^+)^{-1}(p^+(\omega))$, 
		
		\[\supp \mu^+_{p^+(T^{-n}\omega')} \subseteq V_\varepsilon . \]
		Indeed, otherwise, 
		\[\chi(n,\omega')\mu^+_{p^+(T^{-n}\omega')} (U) > \mu^+_{T^{-n}\omega'} (\bdg X  - \overline{V_\varepsilon}) > 0.\]
		Therefore, for all $n$, and for all $\omega' \in \Omega_n \cap (p^+)^{-1}(p^+(\omega))$, $\supp \mu^+_{p^+(T^{-n}\omega')} \cap(\bdg X  - \overline{V_\varepsilon}) = \emptyset$, that is, $\omega' \notin T^n\Omega(V_\varepsilon)$. But $\Omega - T^n\Omega(V_\varepsilon) $ is a subset of measure $\leq \varepsilon$.
		As $\bigcup_n \Omega_n \cap(p^+)^{-1}(p^+(\omega)) $ is conull in $(p^+)^{-1}(p^+(\omega))$, we obtain the set $\omega \in \Omega $ such that $\mu^+_{p^+(\omega)} (U) =0$ is contained in a subset of measure $\leq \varepsilon$. As $\varepsilon$ was arbitrary, we obtain that almost surely, $\mu^+_{p^+(\omega)} (U) > 0$. As this is true for any neighborhood $U$ of $\xi^+$, $\xi^+ \in \supp (\mu^+_\omega)$. Therefore, for all $\omega \in \Omega(V_\varepsilon)$, $\xi^+ \in \supp (\mu^+_\omega)$. As $\varepsilon >0$ was arbitrary, we get the result. 
	\end{proof}

	\begin{cor}\label{cor support limit set}
		For almost every $\omega \in \Omega $, the support of $\mu^+_{p^+(\omega)}$ is the full limit set $\Lambda(G)$. 
	\end{cor}
	\begin{proof}
		The set of attracting points of loxodromic elements is dense in the limit set $\Lambda(G) \subseteq \bdg X$ \cite[Proposition~7.4.6]{das_simons_urbanski17}. By Lemma \ref{lem control support}, the support of $\mu^+_{p^+(\omega)}$ contains almost surely the fixed points at infinity of any loxodromic isometry. Since the support of a measure is closed by definition, we get the result. 
	\end{proof}

	In \cite{maher_tiozzo18}, a great part of the proof for the convergence of the random walk to the boundary amounts to showing that Proposition \ref{prop lim measures hyp} holds. Indeed, the conclusion will follow from the following geometric argument, which we highlight separately. 
	\begin{prop}\label{prop contraction shadow}
		Let $(g_n)$ be a sequence of isometries of an almost geodesic, $\delta$-hyperbolic separable metric space $(X,d)$. Let $\xi \in \bdg X$ and assume that there exists $\xi_1\neq \xi_2 \in \bdg X$ with the following property: for all $U \subseteq \bdg X$, there exists $n_0 $ such that for all $n \geq n_0$, $g_n \xi_1 \in U $ and $g_n \xi_2 \in U $. Then for all $o \in X$, $g_n o \underset{n \to \infty}{\longrightarrow} \xi$. 
	\end{prop}
	\begin{proof}
		Let $S= S_o (x, R)$  be a shadow such that $V:= \bar{S}$ contains $\xi$ in its interior. Let $n_0$ be such that for all $n \geq n_0$, $g_n \xi_1 \in V $ and $g_n \xi_2 \in V $.
		As $X$ is $\alpha$-almost geodesic, Proposition \ref{prop exist qgeod hyp} implies that there exists a $(1,C)$-quasigeodesic line from $\xi_1$ to $\xi_2$, where $C = C(\delta, \alpha)$, which we call $\gamma$. By property of the Gromov product, there is a constant $C = C(\delta, \alpha)$ such that 
		\begin{align}\label{eq dist shadow}
			d(o, \gamma) = (\xi_1, \xi_2)_o + C(\delta, \alpha) =: A.
		\end{align} 
		Now by weak convexity of shadow \ref{cor conv shadows}, there is a constant $C'(\delta, \alpha)$ such that any $(1, C)$ quasigeodesic between points of $S$ is contained in the slightly larger shadow $S_o (x, R + C'(\alpha, \delta))$. Combining this fact with \eqref{eq dist shadow} shows that for all $n \geq n_0$, $g_n o $ is contained in the larger shadow $S_o (x, R+A+ C'(\delta, \alpha))$. As this holds for every shadow $S$ containing $\xi$ in its interior, we have the convergence. 
	\end{proof}

	We shall need the following well-known Portmanteau Lemma, see \cite[Theorem 17.20]{kechris95}. 
	\begin{lem}\label{lem portmanteau}
		Let $Y$ be a separable metrizable space, $P_n$ a sequence of probability measures on $Y$, and $P$ a probability measure on $Y$. Then the following are equivalent: 
		\begin{enumerate}
			\item $P_n \rightarrow_n P$  in the weak-$\ast$ topology. 
			\item $\underset{n \rightarrow \infty}{\liminf} \ P_n(O) \geq P(O)$ for every open set $O \subseteq Y$. 
			\item $\underset{n \rightarrow \infty}{\limsup} \ P_n(F) \leq P(F)$ for every closed set $F \subseteq Y$.
		\end{enumerate}
	\end{lem}
	
	\begin{proof}[Proof of Theorem \ref{thm cv cocycle hyp}]
		Let $\mu^+_\cdot \in \Mapstat_-(\Omega_+, \prob(X))$ be a past stationary measure. Since the action is non-elementary, Lemma \ref{lem control support} implies that the measure $\mu^+_\omega$ is almost surely diffuse. By Proposition \ref{prop lim measures hyp}, for almost every $\omega \in \Omega $, 
		\begin{align}\label{equation qualit contraction}
			\chi(n, \omega)\mu^+_{p^+(T^{-n}\omega)} \underset{n\rightarrow \infty}{\longrightarrow}\delta_{\psi_(\omega)},  
		\end{align}
		in the weak-$\ast$ topology. By Corollary \ref{cor support limit set}, for almost every $\omega$, we have freedom to chose a pair of points $\xi, \eta$ in the support of $\mu^+_{p^+(T^{-n}\omega)}$, disjoint from $\psi(\omega)$. Let $\omega $ be in the intersection of these two conull sets. The qualitative contraction given by Equation \ref{equation qualit contraction} allows us to argue in the same manner as in Proposition \ref{prop contraction shadow}. Let $S= S_o (x, R)$  be a shadow such that $V:= \bar{S}$ contains $\psi(\omega)$ in its interior. Let $\xi_1, \xi_2 \in \bdg X$ disjoint points that are attracting points for two independent loxodromic isometries, and let $S_1, S_2$, be shadows such that $V_1:= \bar{S_1}, V_2:= \bar{S_2}$ contain $\xi_1, \xi_2$ in their interior respectively. As the shadows are disjoint, there exists $C$ such that for all $y_i \in S_i$, $i = 1,2$, $$(y_1, y_2 )_o \leq C.$$ 
		Since $G$ acts by isometries, for all $n$, 
		$$(\chi(n, \omega) y_1, \chi(n, \omega)  y_2 )_{\chi(n, \omega) o} \leq C.$$
		Equation \ref{equation qualit contraction} and the Portmanteau Lemma \ref{lem portmanteau} imply that 
		$$\chi(n, \omega)\mu^+_{p^+(T^{-n}\omega)} (V) \to 1,$$ so that there exists $n \geq n_0$ and $y_1' \in S_1$, $y_2' \in S_2$ such that $\chi (n, \omega) y_1'$ and $ \chi (n, \omega) y_2' $ belong to $ V $. In particular, $\chi (n, \omega) o $ is contained in the larger shadow $S_o (x, R+A+ C'(\delta, \alpha))$. As this holds for every shadow $S$ containing $\psi(\omega)$ in its interior, we have the convergence.
	\end{proof}

\section{Ergodic cocycles in $\cat$(0) spaces}\label{section erg coc cat cv}
	Throughout this section, we let $(\Omega, \mathcal{F}, \mathbb{P}, T, f, G)$ be a RDS on a discrete countable group $G$, and we assume that the measure-preserving map $T$ is invertible and ergodic. We assume that the RDS $\varphi : \mathbb{Z} \times \Omega \rightarrow G$ given by $f$ satisfies the asymptotic past and future independence condition. Let $(X,d)$ be a separable complete $\cat$(0) space and $G \curvearrowright (X, d)$ be a non-elementary action by isometries for which there exists a pair of independent contracting isometries. 
	In particular, the action $G \curvearrowright X$ is non-elementary in the usual sense: there is no $G$-invariant flat (possibly reduced to a point) in $X$. 
	
	Recall that thanks to Proposition \ref{prop non elem},  there exists $L>0$ such that $G$ acts on the hyperbolic model $X_L$ with two independent loxodromic isometries. For the rest of the section, fix such a $L$. We recall that $(X,d_L)$ is $\delta$-Gromov hyperbolic and $\alpha$-almost geodesic \cite[Theorem 3.9]{petyt_spriano_zalloum24}, for constants $\delta, \alpha$ depending only on $L$. To avoid any confusion, in what follows, a geodesic segment $[x,y]$ is always assumed to be in the $\cat$(0) space $(X,d)$, and $B(o, R)$ denotes the open metric ball of center $o$ and radius $R  >0$ for the $\cat$(0) metric $d$. 
	
	\subsection{Convergence to the boundary}
	
	In this section, we prove that the ergodic cocycle $(\chi(n, \omega) o )$ converges almost surely to the visual boundary $\bd X$. This result will follow from Theorem \ref{thm cv cocycle hyp}, as soon as we link the Gromov boundary $\bdg X_L$ of $X_L$ with the visual boundary $\bd X$. We refer to Section \ref{section hyperbolic models} for the vocabulary about the hyperbolic models $(X, d_L)$. 
	
	In the following, even when not specified, the $\cat$(0) spaces we consider are always assumed separable, so that the visual boundary is a standard Borel set.
	
	\begin{thm}\label{thm cv cat}
		Let $(X,d)$ be a separable Hadamard space and $G \curvearrowright (X, d)$ be an action by isometries and assume that $G$ contains a pair of independent contracting isometries for this action. Then, for any $o\in X$, the backward cocycle $(\chi(n, \omega) o)_n$ converges $\mathbb{P}$-almost surely to a point $\xi^+(\omega) \in \bd X$ in the visual boundary. Moreover, there exists $L \geq 0 $ such that almost surely, $\xi^+(\omega) \in B_L$. 
	\end{thm}
	
	\begin{proof}
		Due to Proposition \ref{prop non elem}, there exists $L$ such that $G$ acts on the hyperbolic model $X_L$ by isometries, with a pair of independent loxodromic elements. By Theorem \ref{thm cv cocycle hyp}, the backward cocycle $(\chi(n, \omega) o )_n$ converges almost surely in $X_L$ to a point of the boundary $\xi_L^+(\omega)\in \bdg X_L$. Now due to Proposition \ref{prop inverse map}, $(\chi(n, \omega) o )$ converges almost surely to a point of the boundary $\xi^+(\omega) = \partial_L^{-1}(\xi_L^+(\omega))\in \bd X$. 
	\end{proof}
	\begin{rem}\label{rem separable}
		We stated Theorem~\ref{thm cv cocycle hyp} in the context of roughly geodesic separable hyperbolic spaces. The curtain models are not necessarily separable, but we can pass to an action on a separable hyperbolic space with the arguments of Section \ref{section separable}.
	\end{rem}

	\begin{cor}
		There is a unique invariant measure $\mu_\cdot \in \prob_\bbP(\Omega \times \bd X)$ for the dynamical system $(\Omega, \mathcal{F}, \mathbb{P}, T, f, G)$, which is given by the disintegration $\omega\mapsto \delta_{\xi^+(\omega)}$. Moreover, there is a unique forward-invariant measure $\mu^+_\cdot $. 
	\end{cor}
	
	\begin{proof}
		The first assertion follows from the convergence given by Theorem \ref{thm cv cat}, while the second follows from the one-to-one correspondence between two-sided times and one-sided times in Theorem \ref{thm lim measures rds}.
	\end{proof}
	\subsection{On some contracting properties of the limit points}
	We end this section with a geometric result, that shows that the backward cocycle converges to contracting directions: in case $X$ is proper, limit points are visibility points. 
	
	\begin{prop}\label{prop angle B_L}
		Let $\xi \in \bd X$ be a point of the subspace $B_L$, and let $\eta \neq \xi $ be any other point in $\bd X$. Then $\angle (\xi, \eta) = \pi$. 
	\end{prop}
	
	\begin{proof}
		Let $o \in X$, and let $\gamma : [0, \infty) \rightarrow X$ be a geodesic ray based at $o $ representing $\xi$. By Lemma \ref{lem dual chain}, there is an infinite $L$-chain $\{h_i\}_{i \in \mathbb{N}}$ dual to $\gamma$. Let $\gamma' : [0, \infty) \rightarrow X$ be a geodesic ray representing $\eta$ and such that $\gamma'(0) = o $. By \cite[Remark~8.3]{petyt_spriano_zalloum24}, if $\gamma'$ meets an infinite number of curtains in $\{h_i\}$, it must cross the chain, and hence $\xi = \eta$. As a consequence, there exists $k \geq 1$ such that $\gamma' \subseteq h_{k-1}^{-}$. Let $t_0 >0$ such that for all $t \geq t_0$, $\gamma(t) \in h_{k+1}$. By Lemma \ref{lem bottleneck}, for all $t \geq t_0$ and $t' \in [0, \infty)$, there exists a point $p_{t,t'} \in [\gamma(t), \gamma(t')] \cap h_{k}$ such that, if $q_{t,t'} := \pi_\gamma(p_{t,t'})$, we have $d(p_{t,t'}, q_{t,t'}) \leq 2L+1$. 
		
		Now recall that by \cite[Proposition II.9.8]{bridson_haefliger99}, 
		\begin{align}
			\angle (\xi, \eta) = \lim_{t,t' \rightarrow \infty} \overline{\angle}_o (\gamma(t), \gamma'(t'))  \nonumber. 
		\end{align}
		
		Denote $\overline{\angle}_o (\gamma(t), \gamma'(t')) $ by $\alpha_{t,t'}$. Let $\overline{\Delta} := \overline{\Delta} (\overline{o}, \overline{\gamma(t)}, \overline{\gamma'(t')})$ be the Euclidean comparison triangle of $\Delta (o, \gamma(t), \gamma'(t'))$. By the Euclidean law of cosines, 
		\begin{align}\label{eq al kashi}
			d(\gamma(t), \gamma'(t'))^2 &= d(o, \gamma(t))^2 +d(o, \gamma'(t'))^2 - 2 d(o, \gamma(t))d(o, \gamma'(t')) \cos(\alpha_{t,t'}) \nonumber \\
			&  =  t^2 + t'^2 - 2t t' \cos(\alpha_{t,t'})
		\end{align}
		Since $q_{t,t'}$ belongs to the pole of $h_k$ there exists $A \geq  0$ such that $d(q_{t,t'}, o ) \in  [A, A+1]$ for all $t \geq t_0, t' \geq 0$. By the ``bottleneck'' Lemma \ref{lem bottleneck}, 
		\begin{align}
			t' - (A + 2L + 1) \leq d(\gamma'(t'), p_{t,t'}) \leq t' + A +1 + 2 L+1 \nonumber, 
		\end{align}
		and 
		\begin{align}
			t-( A +1 + 2L + 1 )\leq d(\gamma(t), p_{t,t'}) \leq t+ A+ 2 L+1 \nonumber
		\end{align}
		As a consequence, 
		\begin{align}
			t +t' - (2A+ 4L + 3) \leq d(\gamma'(t'), \gamma(t)) \leq t +t' + 2A + 4 L+3 \nonumber, 
		\end{align}
		Let $C := 2A + 4L +3$. For $t,t' $ large enough, we then have 
		\begin{align}
			t^2 + t'^2 + C^2 + 2 tt' - 2C(t + t') &\leq d(\gamma'(t'), \gamma(t))^2  \nonumber \\
			&\leq t^2 + t'^2 + C^2 + 2 tt' + 2C(t + t')  \nonumber. 
		\end{align}
		Combining this with \eqref{eq al kashi}, we get
		\begin{align}
			C^2 + 2 tt' - 2C(t+t') \leq - 2 tt' \cos(\alpha_{t,t'})
			\leq  C^2 + 2 tt' + 2C(t+t')  \nonumber. 
		\end{align}
		
		As $t, t' \rightarrow \infty$, we obtain that $\cos(\alpha_{t,t'}) \rightarrow -1$, hence $\alpha_{t,t'} \rightarrow \pi$, which proves the proposition. 
		
	\end{proof}
	
	In fact, if the space $X$ is proper, we can show that any point $\xi \in B_L$ is a visibility point.  
	
	\begin{cor}\label{cor angle lim pts}
		Let $\nu \in \prob(\bd X)$ be the hitting measure of the ergodic cocycle $(\chi(n, \omega) o)$. Then for $\nu \otimes \nu$-almost every pair of distinct points $\xi \neq \eta \in \bd X$, $\angle (\xi, \eta) = \pi$. If moreover $X $ is proper, then $\nu$-almost every point of the boundary $\xi \in \bd X$ is a visibility point.
	\end{cor}
	
	\begin{proof}
		The first result is just a combination of Theorem \ref{thm cv cat} and Proposition \ref{prop angle B_L}. The second assertion uses the fact that if moreover the space is proper, points of $B_L$ are visibility points \cite[Lemma 8.4]{petyt_spriano_zalloum24}.
	\end{proof}

	\section{Positive drift and consequences}\label{section drift csq}
	In this section, we study the ``speed'' at which the backward cocycle $(\chi(n, \omega)o)_n$ converges the boundary in both cases. 
	
	\subsection{In hyperbolic spaces}
	
	We let $(\Omega, \mathcal{F}, \mathbb{P}, T,\varphi, G)$ be a greg on a discrete countable group $G$, and we assume that the measure-preserving map $T$ is ergodic. We assume that the RDS $\varphi : \mathbb{Z} \times \Omega \rightarrow G$ given by $f$ satisfies the asymptotic past and future independence condition. Let $(X,d)$ be a separable hyperbolic space and $G \curvearrowright (X, d)$ be a non-elementary action by isometries. 
	We assume here that the cocycle $\varphi$ is integrable, that is, for any $o \in X$, $\omega \mapsto d(\varphi(1, \omega)o, o)$ is integrable: 
	$$\int_\Omega d(f(\omega)o, o) d\mathbb{P}(\omega) < \infty.$$ 
	Observe that for $n,m \in \mathbb{N}$, $\omega \in \Omega$
	\begin{align}
		d(\chi(n+m, \omega)o, o) &\leq  d(\chi(n+m, \omega)o, \chi(n, \omega)o) + d(\chi(n, \omega)o, o) \nonumber \\
		& \leq d(\chi(n, \omega)o, o) + d(\chi(m, T^{-n} \omega) o, o) \nonumber
	\end{align}
	because the cocycle $\chi(n, \omega)$ is backward. In other words, the displacement variables 
	$(D(n, \omega))_n$ defined by
	$$ D(n, \omega) = d(\chi(n, \omega)o , o)$$ 
	form a subadditive cocycle with respect to the flow $T^{-1}$ on $\Omega$. 
	Define the asymptotic average 
	$$\lambda := \lambda_X(\varphi) := \inf \frac{1}{n} \int_\Omega D(n, \omega) d\mathbb{P}(\omega) \in [0, +\infty).$$
	By Kingman subadditive Theorem \cite{kingman68}, the variables $(D(n, \omega)/n)_n$ converge almost surely to $\lambda$ and the convergence also holds in $L^1(\mathbb{P})$. 
	
	\begin{thm}\label{thm drift hyp}
		Assume that the cocycle $\varphi$ is integrable. Then almost surely, $\lambda_X(\varphi) > 0$. 
	\end{thm}
	
	The strategy of the proof is classical. Consider the function 
	\begin{align}
		H\colon & \Omega \times \bdg X \to \bbR \nonumber \\
		 & (\omega, \xi) \mapsto \beta_\xi (f(T^{-1}\omega) o), \nonumber		
	\end{align}
	where $\beta_\xi  = \beta_\xi^o$ is the Busemann function centered at $\xi$ for the basepoint $o \in X$. Define the transformation 
	\begin{align}
		\mathcal{T}\colon & \Omega \times \bdg X \to \Omega \times \bdg X \nonumber \\
		 & (\omega, \xi) \mapsto (T^{-1}\omega, f(T^{-1}\omega)^{-1}\xi). \nonumber		
	\end{align}

	\begin{prop}\label{prop T erg pmp}
		Endow $\Omega \times \bdg X $ with a past invariant measure $\check{\nu} \in \prob_\bbP(\Omega \times \bdg X )$. Then the transformation $\mathcal{T}$ is measure-preserving and ergodic. 
	\end{prop}
	\begin{proof}
		The first statement is by definition. The set of ergodic invariant measures coincides with the set of extreme elements of the convex sets of invariant measures, see \cite[Remark~1.5.6]{arnold98}. Thus ergodicity is a consequence of the classification of invariant measures given by Theorem \ref{thm uniq bdry map hyp}. 
	\end{proof}
	
	Now observe that for all $g_1, g_2 \in G$ and $\xi \in \bdg X$, horofunctions satisfy a cocycle relation: 
	\begin{align}
		\beta_\xi (g_1g_2 o ) & =   \lim_{x_n \rightarrow \xi}d(g_1g_2 , x_n) - d(x_n, o)  \nonumber \\
		& =  \lim_{x_n \rightarrow \xi}d(g_2 , g_1^{-1}x_n) - d(g_1 o, x_n ) +   d(g_1 o, x_n ) - d(x_n, o) \nonumber \\ 
		& =  \lim_{x_n \rightarrow \xi}d(g_2o , g_1^{-1}x_n) - d(o, g_1^{-1}x_n ) +   d(g_1 o, x_n ) - d(x_n, o) \nonumber \\ 
		& =  \beta_{g_1^{-1}\xi} (g_2 o) + \beta_\xi (g_1 o). \label{cocycle horof}
	\end{align}
	Relation \eqref{cocycle horof} gives that for all $\xi \in \bdg X $ and almost every $\omega \in \Omega$,
	
	\begin{align}
		H(\omega, \xi ) + H(\mathcal{T}(\omega, \xi) ) &= \beta_\xi ( f(T^{-1}(\omega))o ) + \beta_{f(T^{-1}(\omega)^{-1}\xi} ( f(T^{-2}(\omega))o ) \nonumber \\
		&= \beta_\xi (f(T^{-1}(\omega))f(T^{-2}(\omega))o) = \beta_\xi (\chi(2, \omega)o),\nonumber 
	\end{align}
	and more generally, that
	\begin{align}
		\beta_\xi (\chi(n, \omega)o) = \sum_{k=1}^{n} H(\mathcal{T}^k (\omega, \xi)) \label{transient cocycle}. 
	\end{align}
	The main geometric ingredient is then the following proposition, which shows that the cocycle $\beta_\xi (\chi(n, \omega)o) $ is a good approximation of the displacement variable $D(n, \omega) $. 
	\begin{lem}
		Let $o \in X$ be a basepoint, and let $(x_n)$ be a sequence converging to $\xi \in \bdg X$ in the Gromov boundary. Then for any $\eta \in \bdg X$ different from $\xi$, and any $o \in X$, there exists $C > 0 $ such that for all $n \geq 0$ we have 
		\begin{equation}
			|\beta_\eta (x_n)  - d(x_n , o)| < C, \nonumber
		\end{equation}
		where $\beta_\eta$ is any horofunction centered on $\eta$. 
	\end{lem}
	\begin{proof}
		Let $(y_p )_p$ be a sequence converging to $\eta$ in the Gromov bordification topology. Observe that 
		\begin{align}
			\beta^o_\eta (x_n) - d(o , x_n) &= \lim_{p \to \infty} d(y_p, x_n) - d(o, y_p) - d(o, x_n ) \nonumber\\ 
			&= \lim_{p \to \infty} -2 (y_p | x_n)_o. \nonumber
		\end{align}
		However, $((y_p | x_n)_o)$ is bounded in $n$ and $ p$ since $\eta\neq \xi$.  
		Then there exists $C >0 $ such that for all $n \in \mathbb{N}$, 
		\begin{equation}
			|\beta^o_\eta (x_n) - d(x_n, o)| \leq C. \label{bdd}
		\end{equation}
		
		Now if we take a different basepoint $o'\in X$, for all $x \in X$, 
		\begin{align}
			|\beta^o_\eta (x_n) - \beta^{o'}_\eta (x_n)| \leq d(o, o'), \nonumber
		\end{align}
		and 
		\begin{align}
			|d(x, o ) - d(x, o')| \leq d(o, o'). \nonumber
		\end{align}
		This shows that the lemma does not depend on the choice of the base-point. 
	\end{proof}
	\begin{cor}\label{prop approx furst cocycle}
		Let $o \in X$ be a basepoint. Then for $\nu$-almost every $\xi  \in \bdg  X$, and $\mathbb{P}$-almost every $\omega \in \Omega$, there exists $C > 0 $ such that for all $n \geq 0$ we have 
		\begin{equation}
			|\beta_\xi (\chi(n, \omega)o)  - D(n, \omega)| < C. \nonumber
		\end{equation}
	\end{cor}

	\begin{proof}[Proof of Theorem \ref{thm drift hyp}]
		Recall that $\beta_\xi $ is 1-Lipschitz on $X$, hence $|H(\omega, \xi)| \leq d(x, \omega_0 x ) $. By the integrability condition on the RDS, $$\int |H(\omega, \xi)| d\mathbb{P}(\omega)d\nui(\xi) < +\infty.$$ 
		
		By Proposition \ref{prop approx furst cocycle}, we have that for $\nu$-almost $\xi \in \overline{X}$, and $\frac{1}{n}\beta_\xi (\chi(n, \omega)o) \rightarrow \lambda$, thus $$\frac{1}{n}\sum_{k=1}^{n} H(\mathcal{T}^k (\omega), \xi) \rightarrow \lambda.$$ In the meantime, due to Proposition \ref{prop T erg pmp}, we can apply Birkhoff Ergodic Theorem \cite[Lemma 2.1]{benoist_quint} and obtain
		\begin{equation*}
			\frac{1}{n}\sum_{k=1}^{n} H(\mathcal{T}^k (\omega, \xi)) \rightarrow \int H(\omega, \xi) d\mathbb{P}(\omega)d\nui (\xi).
		\end{equation*}
		
		Now, by Proposition \ref{prop approx furst cocycle} together with Theorem \ref{thm cv cocycle hyp} gives that $(\beta_\xi(\chi(n, \omega)o)) $ tends to $+\infty$ almost surely. By equation \eqref{transient cocycle}, $(\sum_{k=1}^{n} H(\mathcal{T}^k (\omega), \xi))_n$ is a transient cocycle. Now by the divergence of Birkhoff sums \cite[Lemma 2.18]{benoist_quint}, we obtain that $\int H(\omega, \xi) d\mathbb{P}(\omega)d\nui (\xi) > 0$ so that 
		\[\frac{1}{n}\beta_\xi (\chi(n, \omega)o) \rightarrow \lambda>0.\]
		This finishes the proof of Theorem \ref{thm drift hyp}. 
	\end{proof}
\subsection{In Hadamard spaces}

	\begin{proof}[Proof of Theorem \ref{thm drift cat intro}]
		The drift is positive for apafic ergodic cocycles in hyperbolic models. The result is then just a consequence of Proposition \ref{prop non elem}, and the fact that for any $L\geq o$, the distance in the hyperbolic model $X_L$ satisfies $d_L \leq d$. 
	\end{proof}
	\subsection{Sublinear tracking and sublinear Morse boundary}
	The following is \cite[Theorem 2.1]{karlsson_margulis99}, now that we know that the drift is positive. It states that we have a sublinear geodesic tracking of the ergodic cocycle. 
	
	\begin{cor}\label{cor sublinear tracking}
		Keep the assumptions of Theorem \ref{thm drift cat intro}. Then for $\mathbb{P}$-almost every $\omega \in \Omega$, there is a (unique) geodesic ray $\gamma^\omega : [0, \infty) \rightarrow X$ starting at $o$ such that 
		\begin{equation}
			\lim_{n\rightarrow \infty} \frac{1}{n} d(\gamma^\omega(\lambda n), \chi(n, \omega)o) = 0, \nonumber
		\end{equation}
		where $\lambda $ is the (positive) drift of the ergodic cocycle. 
	\end{cor}
	
	This finishes the proof of Theorem \ref{thm drift cat intro}.

	\bibliographystyle{alpha}
	\bibliography{../bibliography}
\end{document}